\documentclass[12pt]{amsart}

\usepackage{amsfonts,amsmath,latexsym,amssymb,verbatim,amsbsy,amsthm}

\usepackage{graphicx}
\usepackage{nicefrac}
\usepackage{dsfont}
\usepackage{mathrsfs}
\usepackage{cancel}
\usepackage[normalem]{ulem}
\usepackage{color}
\usepackage{pifont}

\usepackage{array}
\usepackage{slashbox}

\usepackage{setspace} 
\usepackage{caption,subcaption} 

\usepackage{algorithm} 

\usepackage{algpseudocode} 

\usepackage{grffile}  

\usepackage[top=1in, bottom=1in, left=1in, right=1in]{geometry}

\usepackage[dvipsnames]{xcolor}
\usepackage{relsize}

\usepackage[colorlinks=true, pdfstartview=FitV, linkcolor=RoyalBlue,citecolor=ForestGreen, urlcolor=blue]{hyperref}

\numberwithin{equation}{section}
\numberwithin{figure}{section}
\numberwithin{table}{section}
\usepackage{float}

\renewcommand{\geq}{\geqslant}

\renewcommand{\leq}{\leqslant}

\newcommand{\ds}{\displaystyle} 
\newcommand{\be}{\begin{equation}}
\newcommand{\ee}{\end{equation}}

\theoremstyle{plain}
\newtheorem{THEOREM}{Theorem}[section]

\newtheorem{theorem}[THEOREM]{Theorem}

\newtheorem{lemma}[THEOREM]{Lemma}
\newtheorem{proposition}[THEOREM]{Proposition}

\theoremstyle{definition}

\theoremstyle{remark}
\theoremstyle{question}

\newtheorem{remark}[THEOREM]{Remark}

\newcommand{\myr}[1]{{\color{black}{#1}}} 

\renewcommand{\a}{\alpha}

\newcommand{\mymin}[1]{{}^{#1}_{-}}  
\newcommand{\mymax}[1]{{}^{#1}_{+}} 
\newcommand{\ioft}{i(t)}
\newcommand{\imin}{{i^-_n}}
\newcommand{\imax}{{i^+}} 

\newcommand{\relmi}{\widetilde{m}_i} 
\newcommand{\relmin}{\widetilde{m}^{n+1}_i} 
\newcommand{\relminin}{\widetilde{m}^{n+1}_{i^-_n}} 
\newcommand{\relminq}{\psi_q(\relmin)} 

\newcommand{\R}{\mathbb{R}}

\newcommand{\calN}{{\mathcal N}}

\newcommand{\lam}{\lambda}
\newcommand{\bX}{\mathbf{X}}
\newcommand{\bx}{\mathbf{x}}
\newcommand{\bxn}{\bx^n}
\newcommand{\xmin}{x^*}
\newcommand{\bxmin}{\bx^*}
\newcommand{\bXnmin}{\bX^n_-}
\newcommand{\bXnplusmin}{\bX^{n+1}_-}
\newcommand{\bXnmina}{\bX^{n_\alpha}_-}
\newcommand{\bXnplusmina}{\bX^{n_\alpha+1}_-}

\newcommand{\bxin}{\bx_i^n}

\newcommand{\hin}{h^n_i}

\newcommand{\etain}{\eta^n_i}
\newcommand{\etainp}{\phi_p(\etain)} 
\newcommand{\bp}{{\mathbf p}}

\newcommand\bpin{\bp^n_i}
\newcommand{\xB}{\bx_B}
\newcommand{\xBi}{(\xB)_i}
\newcommand{\by}{\mathbf{y}}
\newcommand{\mt}{{\mathbf m}^{N}_{t}}
\newcommand{\rd}{\textnormal{d}}
\newcommand{\dt}{{\rd}t}

\newcommand\EE{{\mathbb E}}

\newcommand{\sig}{\sigma}
\newcommand{\gam}{\gamma}

\newcommand{\clam}{\lambda} 

\newcommand{\ddt}{\frac{\textnormal{d}}{\textnormal{d}t}}
\newcommand{\oneoverd}{} 

\def\MAGD{{GD(BT)} }
\def\SGD{GD($h$) }

\DeclareMathOperator*{\argmin}{argmin}
\DeclareMathOperator*{\argmax}{argmax}

\begin{document}

\title[Swarm-Based Gradient Descent]{Swarm-Based Gradient Descent Method\\ for Non-Convex Optimization}
\author{Jingcheng Lu}
\address{Department of Mathematics,\\University of Maryland, College Park}
\curraddr{School of Mathematics,\\University of Minnesota, Minneapolis}
\email{lu000688@umn.edu}

\author{Eitan Tadmor}
\address{Department of Mathematics and Institute for Physical Science \& Technology\newline \hspace*{0.3cm}   University of Maryland, College Park}
\email{tadmor@umd.edu}

\author{Anil Zenginoglu}
\address{Institute for Physical Science \& Technology, 
 University of  Maryland, College Park}
\email{anil@umd.edu}

\date{April 8, 2024}



\thanks{\textbf{Acknowledgment.} Research was supported in part by ONR grant N00014-2112773. ET is indebted to Michelle Tadmor for her comments  that improved the final form of this work.}

\begin{abstract}
We introduce a new Swarm-Based  Gradient Descent (SBGD) method for  non-convex optimization. The swarm consists of agents, each is identified with a position, $\bx$, and mass, $m$.  The key to  their dynamics is  communication: masses  are being  transferred from agents at high ground to low(-est) ground. At the same time, agents change positions  with step size, $h=h(\bx,m)$,  adjusted  to their relative mass: heavier agents  proceed with small time-steps in the direction of local gradient, while lighter agents take larger time-steps based on a backtracking protocol. Accordingly, the crowd of agents is dynamically divided between `heavier' leaders,  expected to approach local minima, and `lighter' explorers. With their large-step protocol,  explorers are expected to encounter improved position for the swarm; if they do, then they  assume the role of `heavy' swarm  leaders and so on. Convergence analysis and numerical simulations in one-, two-, and 20-dimensional benchmarks demonstrate the effectiveness of SBGD  as a global optimizer.  
\end{abstract}

\maketitle
\setcounter{tocdepth}{1}
\tableofcontents

\section{Introduction}

\noindent
 The classical Gradient Descent (GD) methods for optimization, 
 $ \mathop{\argmin}_{\bx\in \Omega\subset \R^d} F(\bx)$, explore the ambient space  by marching along  the directions   dictated by  local gradients, $\nabla F(\bx)$.
  Once the marching direction is determined, the remaining key aspect is  a choice of step size.\newline
More often than not, however,  GD protocols   get trapped in  basins of attraction of local minima, and therefore are not suitable for global optimization of non-convex functions.

\medskip\noindent
In this work we introduce a  swarm-based gradient descent approach  for global optimization.\newline
Communication between agents of the swarm  plays a key role in dictating their step size.\newline
Here,  the usual ambient space of positions is embedded in $\R^d\times [0,1]$: each agent is characterized by its time-dependent position, $\bx_i(t^n)\in \Omega \subset \R^d$, and its relative weight, $\widetilde{m}_i(t^n)\in [0,1]$.\newline 
An interplay between positions and weights proceeds by  communicating a dynamic mass transition  from high to low, thus our protocol favors agents  positioned on  `lower grounds'.\newline 
Looking ahead, the  time-stepping protocol is then adjusted according to  the distinction between `heavy' agents taking  small time steps, and `light' agents taking large(-r)  time steps.\newline
While heavy agents take smaller time-steps, expecting their convergence toward a local minimum, light agents proceed with larger time-steps, so that they explore larger regions, away from local basins of attraction; they are expected to improve the global position of the swarm.  
In the sequel, those light explorers  are expected to encounter a `better' minimizing ground.\newline
Then, these light explorers  are gradually converted into heavier, global leaders of the  swarm.\newline 
Here, the dynamic distinction between heavy leaders and light explorers enables a simultaneous  approach towards local minimizers, while  keep searching for    even better  \mbox{global minimizers}.

\medskip\noindent
Let us recall other well-known multi-agent optimization algorithms based on  `wisdom of the crowd' --- particle swarm optimization  \cite{PSO,grassi2023pso}, \myr{2-agent simulated annealing \cite{chen2019accelerating}}, ant colony optimization \cite{ACO}, genetic algorithms \cite{GA,borghi2023kinetic} and consensus-based optimization  \cite{CBO1,CBO-analytical,carrillo2021consensus,carrillo2023cbo}.\newline
Our Swarm-Based Gradient descent (SBGD) method is shown to be a most effective optimizer, in particular, when the unknown global minimizer is away from the initial swarm.
\mbox{Visiting} larger portions of the ambient space, using explorers based on the communication in swarm dynamics, proved an essential feature for such optimization of remote \mbox{minimizers}. 
Equally important role is played by the leading agents of the swarm: using the backtracking we prove the  sequence of leaders must converge to a minimizer with a quantified rate.

\medskip\noindent
\mbox{Description} of the SBGD method,  given in Section \ref{sec:SBGD},  highlights the decisive  role of communication; indeed, the SBGD can be viewed as  alignment dynamics towards minimal heading. 
\mbox{A precise} time-stepping protocol based on backtracking line search, is outlined in \mbox{Section  \ref{sec:construct}}. 
Detouring the  general paradigm of our swarm-based optimization, we note in Section \ref{sec:outlook} that our recipe for   dynamically adjusting the weights can be extended to  more general \mbox{protocols}.
In particular,  our implementation of SBGD enforces elimination of `worst' agent at each iteration. This `survival of the fittest' approach can be relaxed, increasing the exploring capabilities at the expense of additional computational time. In Section \ref{sec:convergence} we present convergence and error analysis of SBGD. The time-stepping protocol of backtracking implies that the time sequence of SBGD minimizers has a limit set of one (or more) equi-height minima, and depending on the `flatness' of $F$, expressed in terms of Lojasiewicz bound, there follows convergence rate estimate  of the corresponding polynomial order. Finally, in Sections \ref{sec:results-1D}, \ref{sec:results-2D} and \ref{sec:results-20D} we present a series of numerical experiments, comparing the SBGD with various GD methods in one-, two- and respectively 20-dimensional problems. These include GD methods with time-stepping protocol based on a fixed time-step, backtracking and momentum-based Adam protocol \cite{kingma2017adam}. These single-agent methods were implemented using $N$ agents with randomly distributed positions. Of course, having $N$ such agents exploring the region of interest, is  expected to be ``$N$ times better'' than their single-agent versions. Still, when compared with our $N$-based swarm method, we found superior performance of SBGD. Specifically, the communication-based approach in SBGD avoids local minima traps, providing better performance when the search for global minimum requires exploration away from the initial `blob' of  randomly distributed positions.  
\vfill

\section{The Swarm-Based Gradient Descent (SBGD) algorithm}\label{sec:SBGD}
The SBGD dynamics consists of three main ingredients.

\begin{itemize}
\smallskip
\item[\ding{43}] {\bf Agents}. Each agent is identified by its position, $\bx_i(t)\in \R^d$, and its  mass, $m_i(t)\in (0,1]$. The total mass is kept constant in time,   $\sum_i  m_i(t)=1$.

\smallskip
\item[\ding{43}] {\bf Protocol for time step}. The position of each agent is dynamically  adjusted by taking a time step $h_i$ in the gradient direction, $\nabla F(\bx_i(t))$  
\[
\ddt \bx_i(t) = -h_i \nabla F(\bx_i(t)).
\] 
The time step, $h_i$, depends on the position of the agent  at,  $\bx_i(t)$, \underline{and} on its relative mass, $\relmi(t)$,
\[
 \relmi(t):= \frac{m_i(t)}{m_+(t)}, \quad m_+(t)=\max_i m_i(t).
 \]
 The precise dynamic protocol  for choosing the step size,  based on backtracking,  is outlined below. A key aspect is choosing $h_i$  as a decreasing function of the relative mass, $\widetilde{m}_i$: `heavier' agents  move slower, while `lighter' agents  take larger time steps.
An alternative point of view is to interpret  the $\widetilde{m}_i$'s as the probabilities of agents to identify global minimum: those with mass $m_i(t)\ll m_+(t)$ take large time steps to explore the region of interest, since their probability of  identifying the global minimum at their current position, $\bx_i(t)$,  is low.

\smallskip
\item[\ding{43}] {\bf Communication}. Let $\displaystyle F_{\textnormal{max}}(t)=\max_j F(\bx_j(t))$ and $\displaystyle F_{\textnormal{min}}(t)=\min_j F(\bx_j(t))$ denote the maximal and respectively, minimal heights of the swarm at time $t$. The mass of each agent, $m_i(t)$, is dynamically adjusted  according  to its   \emph{relative height}, $\displaystyle \eta_i(t)$,
\[
\qquad \ \ \left\{\begin{split}
\ddt m_i(t)&=-\phi_p(\eta_i(t))m_i(t), \    i\neq \ioft\\
m_i(t) &= 1-\sum \limits_{\ j\neq \ioft}\!\!m_j(t), \   \ \, i=\ioft:=\argmin_i F(x_i(t)),
\end{split}\right. \quad \eta_i(t): = \frac{F(\bx_i(t))-F_{\textnormal{min}}(t)}{F_{\textnormal{max}}(t)-F_{\textnormal{min}}(t)}.
\]  
Thus, each agent `sheds' a fraction of its mass, $\phi_p(\eta_i(t))\in (0,1]$,  which is transferred to the current \emph{global minimizer} at $\bx_{\ioft}$(here we allow to adjust  the mass transition, $\phi_p(\eta)=\eta^p$,  using a user choice of a fine-tuning parameter $p>0$, with the default choice  $p=1$).
As the global minimizer becomes `heavier', it will be `cautious', taking  smaller time steps while enabling the other,  `lighter' agents, to take  larger time steps. As the lighter agents explore the ambient space with larger time steps, it will increase their likelihood to encounter  a new neighborhood of a global minimum, which in turn may  place one of them as the new heaviest global minimizer and so on. Observe that the larger $p$ is, the more tamed the mass transition of $\phi_p(\eta_i(t))$.
 \end{itemize}
 
 \smallskip\noindent 
{\bf  The discrete time} marching of SBGD is realized by agents positioned at $\bx^{n+1}_i=\bx_i(t^{n+1})$ with masses $m^{n+1}_i=m_i(t^{n+1})$ at discrete time steps $t^{n+1}=t^{n}+\Delta t$. We use the simple forward time discretization with  time step $\Delta t=1$, acting on all non-empty agents, $m_i^n>0$, 
\begin{equation}\label{eq:SBGD}
\left\{\begin{array}{l}
      \left.\begin{array}{lll}
         \ m_i^{n+1} & =m^n_i  -\etainp m^n_i, & i\neq \imin \\ \\
          \ m_\imin^{n+1} & =\displaystyle m_\imin^n +\sum \limits_{i\neq \imin}\etainp m_i^{n}, & 
     \end{array} \right\} \quad \imin := \mathop{\argmin}_{i}F(\bxin)\\
     \\
     \quad   \displaystyle m^{n+1}_+:=\max_i m^{n+1}_i \\
     \quad \bx_i^{n+1} = \bxin-h\big(\bxin,\clam\relminq\big)\nabla F(\bxin), \quad \displaystyle \relmin= \frac{m^{n+1}_i}{m^{n+1}_+}
     
\end{array}\right\} \ \ m_i^n>0.
\end{equation}
Initially, the agents are placed  at random positions, $\{\bx^0_i\}$ with  equi-distributed masses $\{m^0_i=\nicefrac{1}{N}\}$.  At each iteration, masses of agents are  exchanged according to their relative heights, and the positions of agents are  updated in the direction of the local gradient, with time step,
$\hin= h\big(\bxin, \clam\relminq\big)$, depending on these relative mass.  Note that the agent with the worst configuration,  positioned at 
$\bx_+=\argmax F(\bxin)$, is eliminated from the computation; consequently, the  size of the swarm decreases, one agent at a time, until it remains with the one heaviest agent. Our choice for the time-stepping protocol, $h\big(\bx,\clam\psi_q(\widetilde{m})\big)$, is the \emph{backtracking line search} outlined in \S\ref{sec:BT}, which is  weighted by the relative masses, $\relminq$ (again, here we allow  fine-tuning the dependence on the relative mass, $\psi_q(\widetilde{m})=\widetilde{m}^q$, based on a user choice of $q>0$, with the default choice $q=1$).
The backtracking enforces a descent property for the SBGD iterations $\bxin$,
and the parameter,  $\clam\in (0,1)$, dictates how much the descent property  holds in the sense that \eqref{eq:SBGDdescent} \myr{below} is fulfilled. The communication  is designed so that the total mass of the swarm gradually concentrates with the agents  most likely to become the global minimizers, that is, the agents which will most likely to reach the global minimum of the region explored so far by the swarm. Such `heavy' agents are assigned with relatively small step sizes, as they are suspected to be close to 'good' minimizers,  hence their subsequent explorations should be sufficiently cautious. On the other hand, the `lighter' agents should not be trapped in basins of attraction of local minimizers, so they proceed with  larger step sizes, allowing them to explore a larger regions, during which they may encounter `better' minimizers; then they may be gradually converted from `light explorers' into `heavy leaders' and so on.
 
We note the flexibility of  the SBGD communication protocol, depending on  fractional mass transition, $\etainp$, and  the mass-dependent step size,  $h(\bxin, \clam\relminq)$. Their detailed construction is outlined in \S\ref{sec:construct}.
 The algorithm \eqref{eq:SBGD} then forms a family of swarm-based methods, denoted SBGD${}_{pq}$ whenever we want to emphasize its dependence on the  parameters $(p,q)$;  the `vanilla' version corresponding to $(p,q)=(1,1)$  is denoted  simply by SBGD. 

\subsection{Why communication is important} 
\begin{figure}[htb]
\centering
     \includegraphics[scale = 0.4]{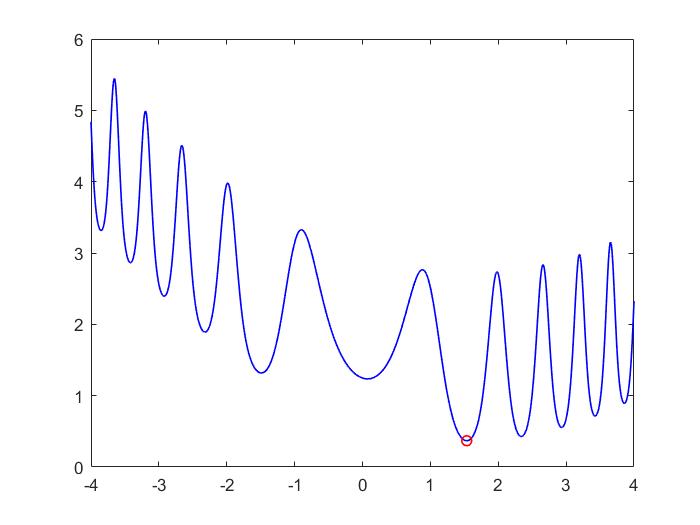}
     \caption{Plot of the objective function \eqref{eq:flat basins}.}
    \label{fig:FB1D}
\end{figure}

Consider the  particular scenario in which all  agents are assigned with the same constant mass, i.e., $m^n_i \equiv \nicefrac{1}{N}$ and $\eta_i^n \equiv 0$,  so that $\psi_q(\widetilde{m}_i^n)\equiv 1$ yields
\begin{equation}\label{eq:GD-BT}
\bx_i^{n+1} = \bxin-h(\bxin,\lambda)\nabla F(\bxin), \qquad i=1,2,\ldots, N.
\end{equation}
 In this case, there is no mass transition  and the dynamics is reduced to a crowd of \emph{non-communicating} agents. In particular, each agent makes its own decision  to proceed with variable step size based on the   backtracking protocol outlined in section \ref{sec:BT} below.
We refer to this  as the backtracking GD, or  \MAGD method. There is also  the vanilla version of  Gradient Descent, denoted GD(h), which proceeds with a fixed step size $h(\bxin,\lambda)\equiv h$. In either case, we have $N$ agents, exploring the region of interest independently of each other. 
Of course, if there are $N$ such agents exploring the region of interest, the corresponding  \SGD and \MAGD method are expected to be ``$N$ times better'' than their the single-agent versions. However,  compared with the  swarm of $N$ communicating agents, we find that  the SBGD dynamics has a superior   \emph{global} behavior. Specifically, the advantage of  communication in SBGD  dynamics becomes  apparent in  exploring  larger regions for potential global minimum. This will be borne out in the numerical results presented in sections \ref{sec:results-1D},\ref{sec:results-2D} and \ref{sec:results-20D}. Here, we  demonstrate  the benefit of communication with a simple example of an  objective function shown in Figure \ref{fig:FB1D},
\begin{equation}\label{eq:flat basins}
    F(x) = e^{\sin(2x^{2})}+\frac{1}{10}\big(x-\frac{\pi}{2}\big)^2.
\end{equation}
The function admits multiple local minima, with a unique global minimum ($x^* \approx 1.5355$).  We compare the performance of SBGD${}_{pq}$,  \eqref{eq:SBGD}  vs. the non-communicating \MAGD iterations, \eqref{eq:GD-BT}, the \SGD iterations with a fixed step size $h$, and the  Adam method with  initial step size $h_0$, denoted Adam($h_0$), \cite{kingma2017adam}. We report on  the results of SBGD${}_{21}$ which  seems to perform slightly better than the `vanilla' version SBGD${}_{11}$, and both offer a more robust optimizer than all other non-communicating methods.

At first, we initialize the positions of agents uniformly in the interval $[-3, 3]$. In this case, the global minimum is included in the support of initial data. We implemented 1000 independent simulations and observe the results of SBGD,  GD(h), \MAGD and Adam. Table \ref{tab:1D FB SBGD vs GD-BT} presents the success rates  for an increasing number of  agents. All methods perform equally well in locating the global minimum, except for Adam(1.1): a large initial step size   in Adam method may take it outside the initial region that already contains the global minimum. 

\begin{figure}[htb]
\centering
    \includegraphics[scale = 0.4]{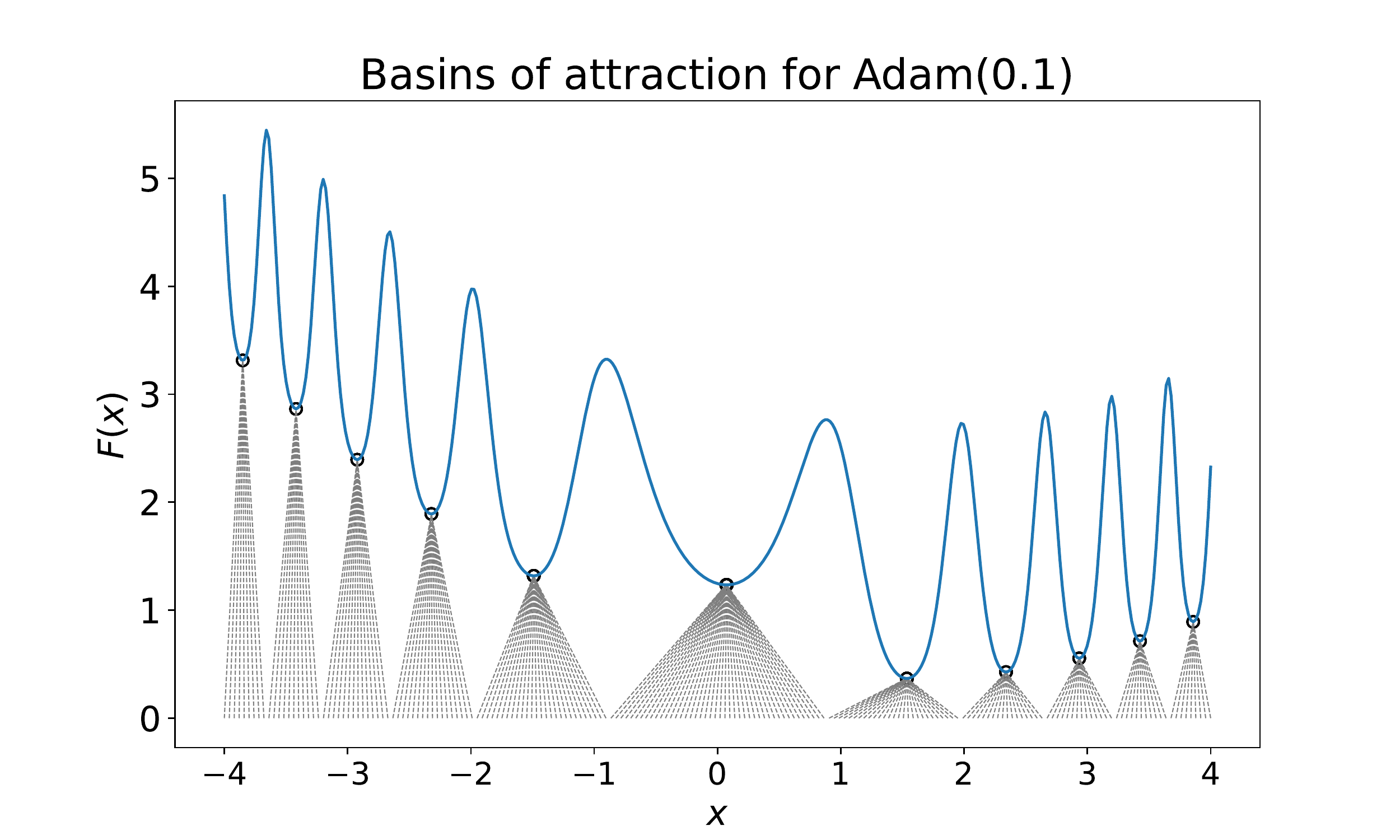}
    \caption{Basins of attraction for GD(0.8) and Adam(0.1) method.}\label{fig:basin of Adam}
\end{figure}

The situation is different, however, if  we initialize    the agents to be uniformly  distributed  in $[-3, -1]$. The results  shown in Table \ref{tab:shifted 1D FB SBGD vs GD-BT}, indicate that the performance of the non-communicating GD(h), \MAGD and Adam(0.1) is significantly worse, whereas the SBGD still identifies the global minimum with high success rates. In particular, the \SGD and Adam with small time steps are trapped inside an initial basin of attraction, unable to get out of that  neighborhood of local minimum. This is also depicted in Figure \ref{fig:basin of Adam} where each local minimum sheds its local basin of attraction for GD(0.8) and Adam(0.1). In particular, the initial data outside $[1,2]$ will necessarily fail to reach the global minimum at $x^* \approx 1.5355$. Only when combined with a larger initial step, Adam(1.1) leads to substantial improvement.

\begin{figure}[h!]
    \centering
    \begin{subfigure}{0.4\textwidth}
    \includegraphics[scale = 0.22]
    {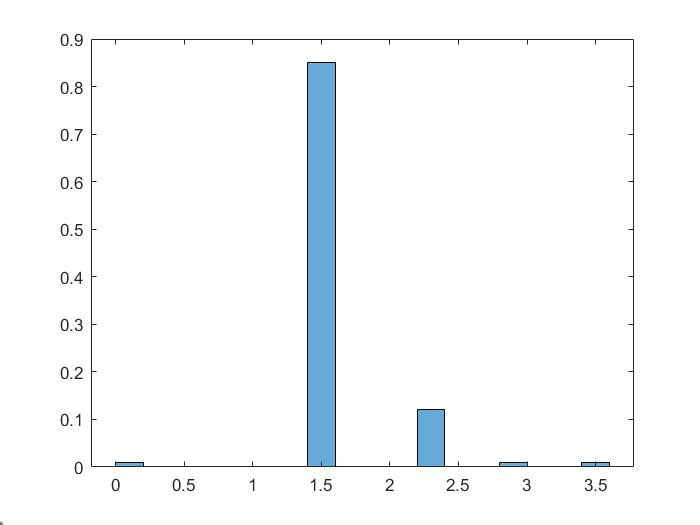}
    \subcaption{SBGD, $N = 10$.}
    \end{subfigure}
    \quad
    \begin{subfigure}{0.4\textwidth}
    \includegraphics[scale = 0.22]
    {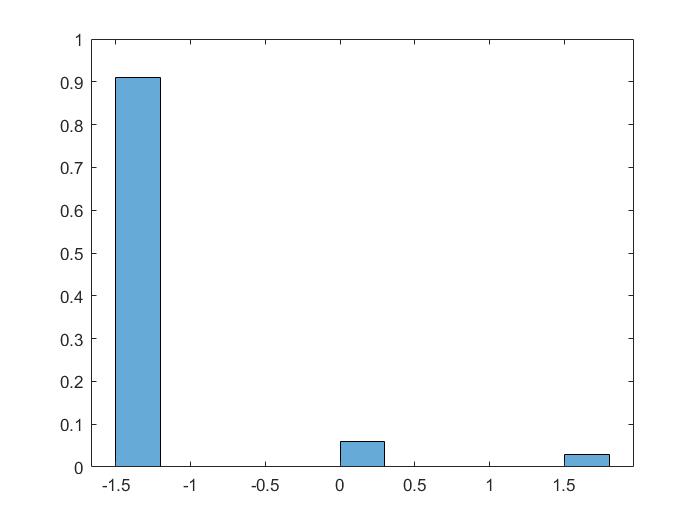}
    \subcaption{\MAGD, $N = 10$.}
    \end{subfigure}
    \\
    \begin{subfigure}{0.4\textwidth}
    \includegraphics[scale = 0.24]
    {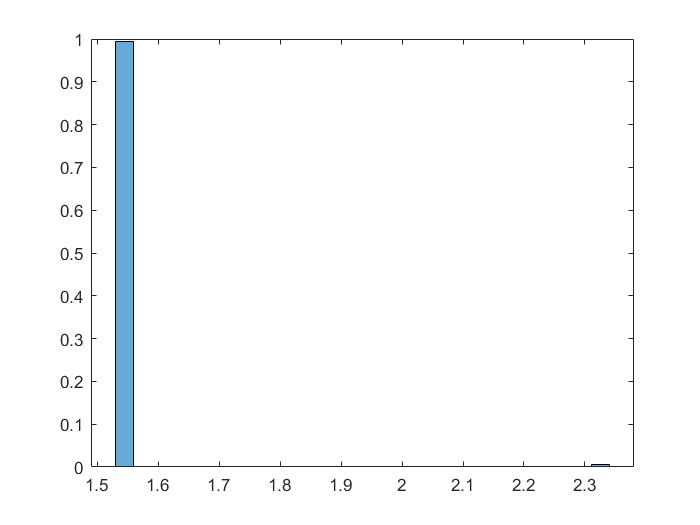}
    \subcaption{SBGD, $N = 20$.}
    \end{subfigure}
    \quad
    \begin{subfigure}{0.4\textwidth}
    \includegraphics[scale = 0.24]
    {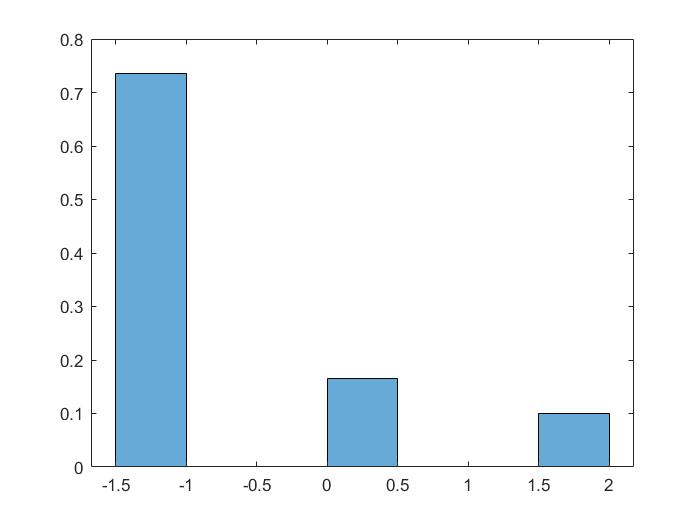}
    \subcaption{\MAGD, $N = 20$.}
    \end{subfigure}
    \\
    \begin{subfigure}{0.4\textwidth}
    \includegraphics[scale = 0.24]
    {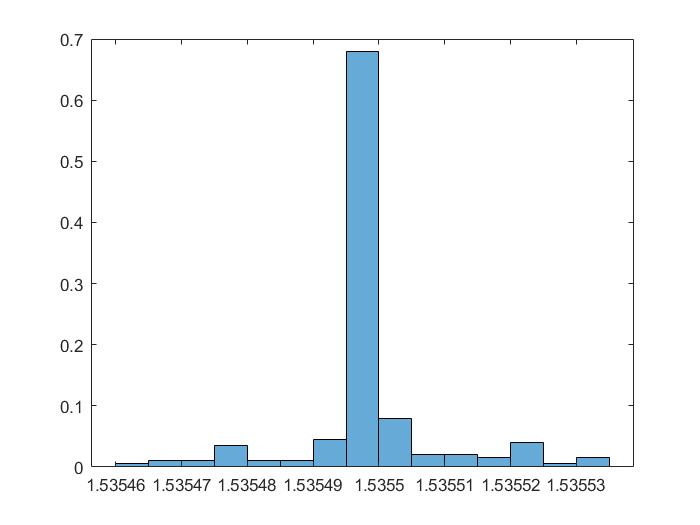}
    \subcaption{SBGD, $N = 30$.}
    \end{subfigure}
    \quad
    \begin{subfigure}{0.4\textwidth}
    \includegraphics[scale = 0.24]
    {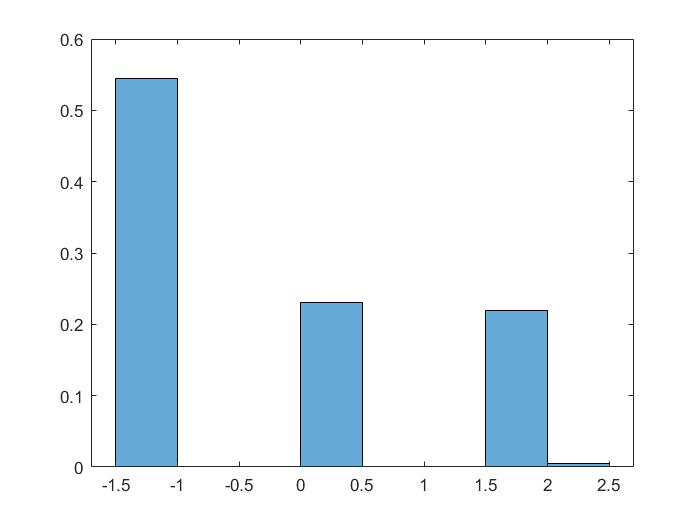}
    \subcaption{\MAGD, $N = 30$.}
    \end{subfigure}
         \caption{Histograms of problem \eqref{eq:flat basins} by $m=200$ experiments. Initial data is generated uniformly in $[-3, -1]$. Global minimum at $\xmin = 1.5355$.}\label{fig:shifted FB1D histogram}
\end{figure}

This is further clarified  when we examine the distribution of $m=1000$ solutions by SBGD vs.  \MAGD   in Figure \ref{fig:shifted FB1D histogram}.    Observe that  in most of the 1000 experiments, the iterations of  \MAGD   are blocked by the  relatively flat basin near the origin, and subsequently they end at the local minimizer lying in the interval $[-2, -1]$. In contrast, the SBGD iterations, thanks to the `aggressive' exploration of light agents, are much more likely to avoid getting trapped in the local flat basin of attraction and eventually accumulate enough mass nearby  the global minimum. 

\begin{table}[htb]
\setstretch{1.5}
       \centering
    \begin{tabular}{|c||c|c|c|c|c|}
    \hline
    N&  5 & 10 & 15 & 20 & 30 \\
    \hline
           SBGD${}_{11}$ & 64.3\% & 96.5\% & 99.8\% & 99.9\% & 100\% \\
     \hline 
     SBGD${}_{21}$ & 68.2\% & 97.7\% & 99.7\% & 100\% & 100\% \\
     \hline 
       GD(0.8)  &  75.2\% & 93.5\% & 98.7\% & 100\% & 100\% \\
       \hline
       \MAGD  &  73.6\% & 96.7\% & 99.5\% & 100\% & 100\% \\
       \hline
       Adam(1.1)  &  19.6\% & 35.0\% & 64.9\% & 77.1\% & 89\% \\
       \hline
       Adam(0.1)  &  58.3\% & 65.6\% & 85.8\% & 95.2\% & 95.7\% \\
       \hline
    \end{tabular}
         \smallskip
 \caption{Success rates  of SBGD, GD(h), \MAGD and Adam methods for global optimization of \eqref{eq:flat basins}, based on   $m=1000$ runs with uniformly generated initial data  in $[-3, 3]$. Backtracking parameters  (see algorithm \ref{alg:backtracking}), $\clam=0.2$ and $\gamma=0.9$.}\label{tab:1D FB SBGD vs GD-BT}
 \end{table}
 
\begin{table}[htb]
\setstretch{1.5}
       \centering
    \begin{tabular}{|c||c|c|c|c|c|}
    \hline
    N&  5 & 10 & 15 & 20 & 30 \\
    \hline
      SBGD${}_{11}$ & 36.5\% & 83.1\% &  97.2\% & 99.5\% &100\% \\
     \hline 
     SBGD${}_{21}$ & 42.4\% & 91.4\% & 99.0\% & 99.8\% & 100\% \\
     \hline 
       GD(0.8) &  0.0\% & 0.0\% & 0.0\% & 0.0\% & 0.0\% \\
       \hline
       \MAGD  &  1.8\% & 5.2\% & 8.5\% & 12.8\% & 21.8\% \\
       \hline
              Adam(1.1)  &  40.2\% & 47.9\% & 82.7\% & 88.7\% & 93.9\% \\
       \hline
       Adam(0.1)  &  0.0\% & 0.0\% & 0.0\% & 0.0\% & 0.0\% \\
       \hline
    \end{tabular}
         \smallskip
 \caption{Success rates of  SBGD, GD(h), \MAGD and Adam methods for global optimization of
  \eqref{eq:flat basins}  based on $m=1000$ runs of uniformly generated initial data in $[-3, -1]$. Backtracking parameters  (see algorithm \ref{alg:backtracking}), $\clam=0.2$ and $\gamma=0.9$}\label{tab:shifted 1D FB SBGD vs GD-BT}
 \end{table}
 
The example shows the advantage of communication-based SBGD as an  algorithm that  is more resilient to the initial guess. Specifically, in complicated applications it may not be realistic to `guess' an initial configuration that encloses the unknown location of the global minimum, and consequently,  GD(h),  \MAGD and Adam iterations may be  trapped near a local minimizer dictated by  ill-conceived initial guesses.  In contrast, the final outcome of SBGD is more resilient with respect to the  initial configuration, in exploring regions outside the enclosure of initial guesses. More can be found in numerical simulations recorded for one-, two- and  20-dimensional benchmark problems presented, respectively, in sections \ref{sec:results-1D}, \ref{sec:results-2D}
and  \ref{sec:results-20D}.  

\subsection{Alignment towards minimal heading} The dynamic adjustment of masses  in SBGD leads to a gradual distinction between `leaders' and `explorers',  according to their relative masses. This adjustment of masses (or probabilities),  $\widetilde{m}_i^n$,  is dictated by  the communication among agents. The SBGD method \eqref{eq:SBGD} can be also interpreted as a particular case of \emph{alignment dynamics}, e.g., \cite{tadmor2021mathematics}, in which agents steer towards the   \emph{minimal} heading, instead of steering towards the average heading \cite{reynolds1987flocks}. In the context of alignment for opinion dynamics, for example, the parameter $\widetilde{m}_i^n$ can be viewed as fraction of the population supporting `opinion' $\bxin$.\newline
This is reminiscent of the Consensus-Based Optimization (CBO) method, 
first proposed in \cite{CBO1} and further modified and analyzed in \cite{CBO-analytical, carrillo2021consensus, CBO-semidiscrete, CBO-timediscrete}; see recent survey \cite{totzeck2022trends}. The CBO method  lets a swarm of agents evolve their positions, $\{\bx^j_t\}$, by a stochastic motion in search of a global minimizer,
 \begin{equation}\label{eq:CBO particle}
    \rd\bx^j_{t} = -\lam(\bx^{j}_{t}-\mt)\dt+\sig(\bx^{j}_{t}-\mt){\rd}{\bf W}^{j}_{t}, \quad \lam,\sig >0
\end{equation}
Agents are driven by two types of motions: the drift towards an exponentially weighted average, $\mt:= \big(\sum_{j=1}^{N}(F_e)^j_t\big)^{-1}\sum_{j=1}^{N}\bx^{j}_{t}(F_e)^{j}_{t}, \ (F_e)^j_t:=\myr{\textnormal{exp}}(-\a F(\bx^{j}_{t}))$, and the stochastic  diffusion, ${\rd}{\bf W}^{j}_{t}$ (implemented by independent Brownian motion different components,  $\sum^{d}_{k=1}(\bx^{j}_{t}-\mt)_{k}{\rd}W^{j,k}_{t}\vec{e}_{k}$). By the Laplace principle, \cite{LaplacePrinciple}, the exponentially weighted average,  $\mt$, concentrates most of its weight with agents of minimal height (or smallest loss). Thus, the drift in \eqref{eq:CBO particle} \emph{aligns} towards those agents with minimal heading, while the stochastic diffusion is responsible for enhancing  the other agents to explore a larger portion of the domain. This should be compared with the deterministic SBGD method, where explorers are driven by communication with lightweight agents.\newline
Unlike  gradient-based methods, the CBO  has the advantage of avoiding  computation of gradients, which are replaced here by a drift towards the $\a$-weighted average. In actual computation, however, in particular in high-dimensional problems,  the CBO method is sensitive to the application of the $\a$-weighted Laplace principle, requiring $\a\gg1$,  which is likely to significantly damage the quality of the solution. It is also sensitive to the choice of parameters $\lambda$ and $\sigma$, \cite{carrillo2021consensus}.

\section{Implementation of the SBGD${}_{pq}$ algorithm}\label{sec:construct}
The success of the SBGD${}_{pq}$ method  relies heavily  on two main procedures: (i) a properly defined communication protocol which dictates mass transition factors, 
$\{\etainp\}$;  and (ii) an effective strategy for taking step size, $h\big(\bxin,\relminq\big)$, which is adjusted to the position and the relative mass of a given agent.
 In this section, we discuss  the details of these procedures, which are summarized in corresponding pseudo-codes.

\subsection{Communications and mass transition}
Let $\displaystyle \imin=\argmin_i F(\bx^n_i)$ denote the global minimizer at iteration number $n$ or time level  $t^n$. All other agents will shed part (or all) of their mass, $\etainp m^n_i, \ i\neq \imin$, which will be transitioned to the mass of the global minimizer $\displaystyle m^n_\imin \mapsto m^n_\imin +\sum_{i\neq \imin}\etainp m^n_i$.   The fraction of mass loss, $\eta^n_i$ is determined by the \emph{relative height}
\begin{equation}\label{eq:etai}
\etain := \frac{F(\bxin)-F^n_{\textnormal{min}}}{F^n_{\textnormal{max}}-F^n_{\textnormal{min}}+\epsilon}, 
\end{equation}
 where $\displaystyle F^n_{\textnormal{min}}= \min_j F(\bx^n_j)$ and $\displaystyle F^n_{\textnormal{max}}= \max_j F(\bx^n_j)$ are the current global extremes\footnote{To  prevent vanishing denominator in the extreme case $F_{\textnormal{max}} = F_{\textnormal{min}}$, we introduce a small $\epsilon$-correction, say $\epsilon = 10^{-10}$;  the effect on the numerical performance is minimal.}. In this fashion, the `higher' the agent is, the more mass it will lose, and indeed, the highest agent in each iteration will be eliminated\footnote{To be precise,  the worst agent is eliminated whenever  its relative height $\eta^n_i=1$. This can be realized only when $\epsilon=0$.  When $\epsilon >0$, elimination of worst agents takes place whenever $1-\eta^n_i = {\mathcal O}(\epsilon) \ll 1$.}.
 The function $\phi$ is user-dependent; for example $\phi_p(\eta)=\eta^p, \ p>0$ enables adjusting the amount of mass transition, with  mass transition tamed as $ p \uparrow \infty$. 
 This usage of a relative height keeps a minimal amount of  global communication necessary to calibrate  each agent relative to the current extremes of the crowd, while being invariant under the transition and dilation of the target function. Therefore, the computation of the relative height is more stable.

 \subsection{Backtracking -- a protocol for time stepping}\label{sec:BT}
 Consider the vanilla Gradient Descent (GD) iteration
 \[
 \bx^{n+1}=\bx^n -h\nabla F(\bx^n).
  \]
 The new position, $\bx^{n+1}=\bx^{n+1}(h)$, is viewed as a function of the step size, $h$. A proper strategy for choosing the step size  is now the key for the success of GD iterations. We recall the classical  \emph{backtracking line search},   \cite[\S3]{nocedal2006conjugate}, which is a computational realization of the  well-known Wolfe conditions \cite{wolfe1969,armijo1966} for inexact line searches. The idea  is to secure an acceptable step length, $h$, that enforces a sufficient amount of  height reduction (or loss) in the target function  
\begin{equation}\label{eq:backtracking}
  F(\bx^{n+1}(h)) \leq F(\bxn)- \clam h |\nabla F(\bxn)|^2, \quad \clam\in (0, 1).
\end{equation}
Of course, since $F((\bx^{n+1}(h)) = F(\bxin)- h  |\nabla F(\bx^n)|^2+\textnormal{\{higher  order  terms\}}$, then 
\eqref{eq:backtracking} holds for any fixed $0<\lambda <1$, provided the time step is small enough, $h\ll 1$. The purpose is to secure \eqref{eq:backtracking} for  large enough $h$,  so that we maximize the size of descent, $\clam h|\nabla F(\bxn)|^2$.
To this end, one employs a dynamic adjustment, starting with a relatively large  $h$ (for which one expects
$\displaystyle 
F\big(\bxn-h\nabla F(\bxn)\big)>F(\bxn)- \clam h|\nabla F(\bxn)|^2$)
  and then \myr{repeatedly} shrink the step size \myr{$h \rightarrow \gamma h, \ \gamma<1$} until  \eqref{eq:backtracking} is observed.  At this stage, we have a final step size, $h=h(\bxn, \clam)$, which secure \eqref{eq:backtracking}
  \[
  F(\bx^{n+1}(h)) \leq F(\bxn)- \clam h(\bxn,\clam) |\nabla F(\bxn)|^2.
  \]
 Thus, the protocol for step size  hinges on the choice of the   parameter $\clam$ which  dictates the amount of  descent property \eqref{eq:backtracking}. One expects that as the iterations approach the (potentially) global minimum their descent property is `tamed' with a larger $\clam$, yet they  should be able to avoid getting trapped in local basins of attraction by allowing smaller $\clam \ll 1$. This is exactly where we take advantage of our swarm-based approach: as a compromise between these two conflicting requirements, we offer to use the relative mass of different agents, $\relmin$, as an indicator which distinguishes between the `heavier' agents who are potentially  close to  the global minimizer, and the `lighter' agents which are  allowed to take large(r) time steps.  
We therefore adjust the descent parameter $\clam$ to each  SBGD${}_{pq}$ agent positioned at   $\bxin$, according to its relative mass $\relmin$ 
\begin{equation}\label{eq:step length}
    \hin=h\big(\bxin,\clam\relminq\big)  \qquad \relmin:=\frac{m^{n+1}_i}{m\mymax{n+1}}, \quad  m\mymax{n+1}:= \max_{i} m^{n+1}_i, \ \ 0<\clam<1.
\end{equation}
The parameter $q>0$ tunes the role of the relative mass $\relmin$: as $q$ increases then $\relminq$ decreases, and backtracking allows  intermediate agents to take larger time steps.
\noindent
The pseudo-code for computing the SBGD${}_{pq}$ steps based on backtracking line search is given in Algorithm \ref{alg:backtracking}.
The results reported in \S\ref{sec:pq} below show that although fine-tuning the parameter, $q$, can lead to improved results,  it has a limited effect on the overall performance of SBGD${}_{pq}$ iterations.

\begin{algorithm}[h!]
\setstretch{1.5}
\caption{Backtracking Line Search}
\label{alg:backtracking}
 \begin{flushleft}
 \% \myr{Determine $\hin$:  the SBGD step-size of the agent positioned at $\bx_i$ at time $t^n$}
 
Set the descent parameter, $\clam \in (0, 1)$, and shrinkage parameter, $\gam\in(0, 1)$

\noindent
Set $\psi_q$ with $q>0$

\noindent
Set the relative mass $\displaystyle \relmin=\frac{m^{n+1}_i}{m^{n+1}_+}$

\noindent
Initialize the step size, $h = h_0$, \myr{with  large enough $h_0$ (see \eqref{eq:hlower-bound} below}).
\end{flushleft}

\begin{algorithmic}
\While{$F\big(\bxin-h\nabla F(\bxin)\big)>F(\bxin)- \clam\relminq h|\nabla F(\bxin)|^2$}

\State $h\gets \gamma h$.

\EndWhile

 \begin{flushleft}
Set $\myr{\hin} \gets h$ \qquad \qquad \% sets $\hin$ as the step size depending on $\bxin$ and $\clam\relminq\big)$
\end{flushleft}
\end{algorithmic}

\end{algorithm}
\noindent
Observe that the successive shrinking of  step size in backtracking involves a shrinking factor $\gamma \in (0, 1)$: a  small $\gamma$ corresponds to a `crude' line search while $\gamma\sim 1$, corresponds to a more refined search. Again, although  fine-tuning the shrinkage parameter $\gamma$ may lead to improved results, it does not seem to have a substantial effect on the overall performance of SBGD${}_{pq}$. \myr{It does, however, comes at a  substantial cost: a fine-tuning of $\gamma \sim 1$ may end up with  many function calls of $F$-evaluations, before  its descent bound, analyzed in Lemma \ref{lemma:lower bound for backtracking} below,  is fulfilled. On the other hand, such tighter `screening' of the ambient space is more likely to lead to better, i.e., lower minima values.}. 

\medskip
The final output of algorithm \ref{alg:backtracking} yields,  for each agent, an adjusted  step size, $\hin=h\big(\bxin,\clam\relminq\big)$, which secures the descent property \eqref{eq:backtracking} with $\clam\relminq$  substituted for $\clam$,
\begin{equation}\label{eq:SBGDdescent}
F(\bx^{n+1}_i) \leq F(\bxin)- \clam\relminq \hin |\nabla F(\myr{\bxin})|^2, \qquad  \quad \bx^{n+1}_i= \bxin-\hin\nabla F(\bxin).
\end{equation}
Moreover, as we shall see in Lemma \ref{lemma:lower bound for backtracking} below, the step sizes, $\hin$, admit a lower bound in terms of the corresponding relative masses of the different agents.  
The scaling of the step size  using relative masses encodes the communication  of different agents, which is the key to the success of the SBGD${}_{pq}$ algorithm.
Roughly speaking, we can distinguish between two types of agents. The SBGD${}_{pq}$ iterations are led by the heavier agents,  $m^{n+1}_i \approx m{\mymax{n+1}}$, which tend to recover a maximal  local descent  rate of order $\clam$, \eqref{eq:backtracking}. On the other hand, there are the lighter agents where $\displaystyle m^{n+1}_i \ll m{\mymax{n+1}}$, which are less driven by the steepness of their decent, and are therefore better equipped as  explorers  of large areal search for the global minimizer. In this way, the mass-dependent adjustment of step size  captures both the descent property of  the target function while allowing the lighter explorers to pull away from local minimizers.

\subsection{SBGD${}_{pq}$ pseudocode}\label{sec:implement}
The pseudocode of the SBGD${}_{pq}$ method is given in Algorithm \ref{alg:SBGD}. The initial setup consists of $N$ randomly distributed agents $\bx^{0}_{1},\cdots,\bx_{N}^{0}$, associated with masses $m_{1}^{0},\cdots,m_{N}^{0}$. At the beginning, all  agents are assigned  equal masses, $\displaystyle m_{j}^{0} = \nicefrac{1}{N}$, $j = 1,\ldots, N$. At each iteration step, the agent $\bx_\imin =  \argmin_{\bxin} F(\bxin)$ attains the minimal value, while the other agents  transfer part of their masses to the current optimal minimizer $\bx_\imin$. Then all the agents are updated with the gradient descent method using the step lengths obtained with \eqref{eq:step length}. 

\noindent
To further improve  efficiency, we use three tolerance factors:\newline
$\cdot$ If the mass of an agent is lower than a minimal threshold $tolm$, then this agent will be eliminated and its remaining mass will be transferred to the optimal agent at $\bx_\imin$.\newline 
$\cdot$ ``Sticking particles''. Agents that are sufficiently close to each other below a threshold $tolmerge$, are merged into a new agent, and  their masses are combined into  the newly generated agent.\newline
$\cdot$ The iterations stop when the  minimizer's descent in two consecutive iterations is below a minimal threshold $tolres$.

\noindent
Unless otherwise specified, all simulations  reported in this paper employ the same thresholds
\begin{equation}\label{eq:tols}
tolm = 10^{-4}, \quad tolmerge = 10^{-3}, \quad tolres = 10^{-4}.
\end{equation}
The optimal choices of these thresholds are experimental.  With smaller thresholds, $tolm$ and $tolmerge$,  the SBGD algorithm will explore a larger part of the ambient space ending with a better solution, at the expense of  reduced efficiency. A balance between the quality of the solution and the computational cost should be explored.

\begin{algorithm}[h!]
\setstretch{1.1}
\caption{Swarm-Based Gradient Descent}\label{alg:SBGD}

\begin{algorithmic}

\State Set three tolerance parameters, $tolm$, $tolmerge$ and $tolres$

\State Set the adjustment parameters $p, q>0$

\State Initialization:

\State \qquad Set $N$ --- the number of agents

\State \qquad Set  initial positions $\bx_{1}^{0},\cdots,\bx_{N}^{0}$ randomly generated under initial distribution $\rho_0$

\State \qquad Set initial mass $\displaystyle m^{0}_{1}=\cdots=m^{0}_{N} = \nicefrac{1}{N}$

\State \qquad Set  the optimal agent, $\displaystyle i^-_0 = \argmin_{i} F(\bx^0_i)$

\For{$n = 0,1,2,\cdots$}

\State Set $\displaystyle F{\mymin{n}} = F(\bx^n_\imin)$, $\displaystyle F{\mymax{n}} = \max_{i}F(\bxin)$

\For{$i = 1,\cdots,N$ and $i\neq \imin$} \hspace{0.6em} \% Mass transitions

\If{$\displaystyle m_i^n<\nicefrac{1}{N}*tolm$} 
\State set $m_i^{n+1}= 0$
\State reduce the $\#$ of active agents: $N\gets N-1$
\Else  \ $\displaystyle m_i^{n+1} = m^n_i -\etainp m^n_i$  where 
$\displaystyle \etain = \frac{F(\bxin)-F^n_{\textnormal{min}}}{F^n_{\textnormal{max}}-F^n_{\textnormal{min}}}$.

\EndIf
\EndFor

\State $\displaystyle m_\imin^{n+1} = m_\imin^n +\sum_{i\neq \imin} \etainp m^n_i$ 
\hspace{0.5in} \% The mass of the overall crowd is conserved

\State Compute $\displaystyle m_+=\max_i m^{n+1}_i$
\For{$i = 1,\cdots,N$} \hspace{0.7in} \% Gradient descent

\State Compute relative masses $\displaystyle \relmin=\frac{m^{n+1}_i}{m_+}$

\State Compute the step size $h = h\big(\bxin, \clam\relminq\big)$ according to algorithm \ref{alg:backtracking}.

\State March: $\displaystyle \bx_i^{n+1} = \bxin-h \nabla F(\bxin)$.

\EndFor


\State Merge the agents if their distance $<tolmerge$. 

\State Set  the new optimal agent $i^-_{n+1}$ = $\displaystyle \argmin_{i} F(\bx^{n+1}_i)$.

\State Compute the residual  $\displaystyle res = |\bx^{n+1}_{i_{n+1}}-\bx^n_{i_n}|_{2}$

\If{$res<tolres$} 
\State $\bx_{SOL}  \gets \bx^{n+1}_{i^-_{n+1}}$
\State break
\EndIf

\EndFor
\end{algorithmic}
\end{algorithm}

\section{A general outlook}\label{sec:outlook}
We are aware that there are  many possible extensions that can be worked out in connection with the SBGD algorithm, leading to a large class  of swarm-based optimizers (SBO) with better communication protocols.
We mention three of them.

\smallskip
\paragraph{$\cdot$ \bf General gradient descent directions}
Our SBGD approach can be used  with  a  more general set of gradient descent directions, $\bpin$

\[
\bx^{n+1}_i=\bxin + h\big(\bxin,\clam\relminq\big) \bpin, \qquad \langle \bpin, \nabla F(\bxin)\rangle <0.
\]
The convergence results in \S\ref{sec:convergence}, with appropriate adjustments, remain valid. Moreover, the variety in choice of directions, other than local gradients, may offer a better `covering' of the ambient space 
$\Omega \subset {\mathbb R}^d$.  \myr{This  swarm-based descent approach with \emph{randomly} chosen descent directions satisfying $\langle \bpin, \nabla F(\bxin)\rangle <0$ was pursued in our recent work \cite{tadmor2024swarm}.}

\smallskip
\paragraph{$\cdot$ \bf Swarm-based optimization --- a general paradigm}
The general paradigm for our swarm-based optimization is realized  by embedding the $d$-dimensional ambient space in $(\bx,\widetilde{m})\in {\mathbb R}^d\times [0,1]$; here $\widetilde{m}$ is an additional parameter space of masses/weights (or probabilities,  or `fractional population', ...) which serves as  a communication platform  for the crowd of agents positioned at $\{\bx_i\} \in {\mathbb R}^d$.  In this context, one can combine such communication-based swarm  iterations with any single agent time-marching protocol. As   examples we refer to the recent 
adaptive GD method \cite{liu2022adaptive} and the references therein.  In the present work we use  the time-marching protocol of gradient-descent,  based on backtracking search. Other time marching protocols can be used.

\smallskip
\paragraph{$\cdot$ \bf Survival of the fittest}
The communication in SBGD is designed so that in each iteration, the `worst' agent, positioned at  $\displaystyle \bx_{\imax} := \argmax_{\bx_i}F(\bxin)$, is eliminated,  as it loses all of its mass ($\eta_\imax^{n} = 1$). This  policy of  `survival of the fittest' implies that  the number of $N$ initial active   agents  decreases in each  iteration until the SBGD   remains with  only a single,  `heaviest' agent, which proceeds by the \MAGD protocol. In particular, this policy implies that for small swarms, say $N\sim 10$, the performance of SBGD  iterations is expected to be similar or only slightly better than \MAGD iterations, as borne out in the numerical simulations reported in sections \ref{sec:results-2D} and \ref{sec:results-20D}.  Alternatively, one can design a less restrictive evolutionary policy that will allow `worst' agents below a certain threshold to survive. This will evolve a larger set  of explorers for longer times, with a greater chance of exploring new and better minima unseen before. Our numerical experiments show that a balanced policy for the `fittest' can indeed have a substantial effect on the final result, at the expense of 
increased  computational time. \myr{The main target of these different policies is to computational efficiency; we emphasize that the convergence analysis carried out on \S\ref{sec:convergence} is independent of the `evolutionary policy'.}

\section{Convergence and error analysis}\label{sec:convergence}
The study of convergence and error estimates for the SBGD method requires  to quantify the behavior  of $F$. Here we emphasize  that the required smoothness properties of $F$ are only sought in the  region explored by  the  SBGD iterations. We  assume  that there exists a \emph{bounded} region, $\Omega \ni  \bxin$ for all agents. Since the SBGD allows light agents to explore the ambient space with large step size (starting with $h_0$), we do not have apriori bound on $\Omega$; in particular, the footprint of the SBGD crowd, $\myr{\cup_n}\textnormal{conv}_i\{\bxin\}$, may expand well beyond its initial convex hull
$\textnormal{conv}_i\{\bx^0_i\}$. \myr{We let $\bx^*:=\argmin_{\bx\in \Omega} F(\bx)$ denote the global minimum value in that region. This is the  minimum that we would wish to converge to}.
\newline
We consider the class of loss functions,  $F\in C^2(\Omega)$, with  Lipschitz bound
\begin{equation}\label{eq:FisC2}
    |\nabla F(\bx)-\nabla F(\by)|\leq L|\bx-\by|,\quad \forall \bx,\by \in \Omega.
\end{equation}
We begin by recalling the lower bound on the step size, secured by the 
 backtracking line search in algorithm \ref{alg:backtracking}.
 \begin{lemma}\label{lemma:lower bound for backtracking}
  Consider the SBGD${}_{pq}$ iterations  \eqref{eq:SBGD}, with  step size $\hin=h\big(\bxin,\clam\relminq\big)$ determined by the backtracking line search in algorithm, \ref{alg:backtracking}, with shrinkage factor $\gam\in (0,1)$ and  initial step size, $h_0$, large enough so that
 \begin{equation}\label{eq:hlower-bound}
 \myr{h_0 > \frac{2}{L}}.
 \end{equation}
 Then we have the descent bound
 \begin{equation}\label{eq:allisx}
\myr{F(\bx^{n+1}_i)   \leq F(\bxin)-\frac{2\gam}{L}\big(1-\clam\relmin\big)\clam\relmin|\nabla F(\bxin)|^2, \qquad \bx^{n+1}_i=\bxin -\hin\nabla F(\bxin).}
\end{equation}
\end{lemma}

\noindent
\myr{We make two comments before turning to the proof of the lemma.\newline
First we note that if one employs a step size $h$ small enough then there holds
\begin{equation}\label{eq:if-h-small-enough}
F\big(\bxin -h\nabla F(\bxin)\big)   \leq F(\bxin)-\frac{2\gam}{L}\clam\relmin h|\nabla F(\bxin)|^2.
\end{equation}
However, in this case the `amount of descent', $\ds \frac{2\gam}{L}\clam\relmin h|\nabla F(\bxin)|^2$, is comparably small with $h \ll 1$. In contrast, \eqref{eq:allisx} secures the amount of descent which depends on the relative mass, but otherwise is independent of the smallness of $h$. To this end, we need to secure a step size $\hin$ \emph{large enough}; this is the goal of the backtracking algorithm \ref{alg:backtracking}.\newline
Second --- it is customary to find stability results which state that {``if $h$ is small enough than such and such''}; typically in this case, e.g., \eqref{eq:if-h-small-enough}, ``smallness'' is measured by $1/L$ and the difficulty is that in general we have no access to how large $L$ is. The point of Lemma \ref{lemma:lower bound for backtracking} is that the descant  \eqref{eq:allisx} holds \emph{without} access to $L$ (to be precise,  \eqref{eq:hlower-bound} requires $1/L$; however, here $h_0$ needs to be \emph{large} relative to $1/L$)}. 

 \noindent \emph{Proof of Lemma \ref{lemma:lower bound for backtracking}.} 
We claim that the descent bound holds
\begin{subequations}\label{eqs:back}
\begin{equation}\label{eq:back}
F(\bx^{n+1}_i)\leq F(\bxin)-\clam\relminq\hin|\nabla F(\bxin)|^2, \qquad \
\end{equation}
with step size satisfying the lower bound
\begin{equation}\label{eq:lower bound for hk}
\hin \geq \frac{2\gamma}{L}\big(1-\clam\relminq\big).
\end{equation}
\end{subequations}
Clearly, \eqref{eq:allisx} follows from \eqref{eqs:back}.\newline
Now, it may happen that \eqref{eq:back} holds for our initial choice of a step size $h_0$ which already tuned to satisfy \eqref{eq:lower bound for hk}, in which case we are done. But in general, $h_0$ need not secure any descent at all and we appeal to the backtracking iterations.
By the Lipschitz continuity of $\nabla F$, 
\begin{equation*}
\begin{split}
F\big(\bxin-h\nabla F(\bxin)\big)&\leq F(\bxn)-h|\nabla F(\bxin)|^2+\frac{L}{2}h^2|\nabla F(\bxin)|^{2}   \\
& = F(\bxin)-\Big(1-\frac{L}{2}h\Big)h|\nabla F(\bxin)|^2,
\end{split}    
\end{equation*}
and hence, if $h$ is small enough  
\[
    h \leq \frac{2}{L}\big(1-\clam\relminq\big) \ \ \leadsto \ \ F(\bxin-h\nabla F(\bxin))\leq F(\bxin)-\clam\relminq h|\nabla F(\bxin)|^2.
\]
The backtracking line search iterations tell us that the inequality on the right holds for $h=\hin$, that is \eqref{eq:back} holds;  but it does not hold for $\ds \frac{\hin}{\gamma}$ (the stopping criterion fails with $\hin/\gamma$). In particular, therefore, $\ds \frac{\hin}{\gamma}$ must satisfy the reverse inequality on the left,
that is, \eqref{eq:lower bound for hk} holds. $\square$

\subsection{Convergence to a band of local minima} Our next proposition provides a rather precise quantitative description for the convergence of the SBGD method. The convergence   is determined by the time series of  SBGD minimizers, $\{\bXnmin\}$,
\[
\bXnmin=\bx^n_{i^-_n}, \qquad i^-_n:=\argmin_i F(\bxin),
\]
 and the time series of  \myr{\emph{parents}} of the heaviest agents, $\{\bX\mymax{n}\}$, 
 \[
 \bX\mymax{n}=\bx^n_{i^+_{n+1}}, \qquad i^+_{n+1}:=\argmax_i m^{n+1}_i.
 \]
 \myr{Note that $\bX\mymax{n}$ is the ``parent'' of $\bx^{n+1}_{i^+_{n+1}}$ which is, by definition, the heaviest agent at $t=t^{n+1}$}.
 The interplay between minimizers and communication of masses leads to a gradual   shift of mass, from higher ground to the minimizers. Eventually, when the SBGD minimizers gain enough mass to assume the role of heaviest agents, the two sequences coincide.  
Convergence is independent of the lighter agents.\newline
We  introduce the  scaling $M=\max_j F(\bx_j^0)-F(\bxmin)$; since $F(\bxin)$ are decreasing, we conclude that the SBGD iterations remain within that range, namely
\begin{equation}\label{eq:var}
 \forall n, i: \quad \max_j F(\bx^n_j)-F(\bxin) \leq M, \qquad M:=\max F(\bx_i^0)-F(\bxmin)
 \end{equation}
To simplify matters we restrict our attention to the vanilla version of SBGD, $(p,q)=(1,1)$.
\begin{proposition}\label{prop:convergence of SBGD}
Fix $\clam \in (0,1)$ and consider the SBGD iterations  \eqref{eq:SBGD} with  step size $\hin=h\big(\bxin,\clam\relmin\big)$ determined by backtracking line search of algorithm \ref{alg:backtracking} in \eqref{eq:step length}.\newline
  Let $\{\bXnmin\}_{n\geq0}$ and $\{\bX\mymax{n}\}_{n\geq0}$ denote the time sequence of SBGD  minimizers  and, respectively,   (parent of) heaviest agents, at $t^n$.
 Then, there exists a constant, $C=C(\gamma,L,M, \clam)$ given in \eqref{eq:telesC} below, such that
\begin{equation}\label{eq:summability}
   \myr{ \sum_{n = 0}^{\infty} |\nabla F(\bX\mymax{n})|^{2} \times \min\big\{1,|\nabla F(\bXnmin)|^{2}\big\}  < CM}.
\end{equation}
\end{proposition}

\noindent \emph{Proof.} 
By Lemma \ref{lemma:lower bound for backtracking}, 
the descent property \eqref{eq:allisx}  of different agents is dictated by their relative mass
\begin{equation}\label{eq:allis}
F(\bx^{n+1}_i)  \leq F(\bxin)-\frac{2\gam}{L}\big(1-\clam\relmin\big)\clam\relmin|\nabla F(\bxin)|^2, \qquad \relmin=\frac{m^{n+1}_i}{m^{n+1}_{i^+_{n+1}}}.
\end{equation}
In particular, application with $i=i^+_{n+1}$  implies the descent of the heaviest agent positioned at $\bX\mymax{n}=\bx^{n}_{i^+_{n+1}}$
\begin{equation}\label{eq:heaviest}
F(\bx^{n+1}_{i^+_{n+1}})  \leq F(\bX\mymax{n})-\frac{2\gam}{L}\big(1-\clam\big)\clam|\nabla F(\bX\mymax{n})|^2.
\end{equation}
We also need  to  quantify the descent  of  the minimizer 
 positioned at $\bXnmin=\bx^{n}_{i^-_{n}}$; this requires a lower bound on its relative mass $\displaystyle m^{n+1}_\imin$ on the right of \eqref{eq:allis}.
To this end, we  consider two sub-cases, depending on the size of $\ds F(\bX\mymax{n})-F(\bXnmin)>0$.\newline
Case (i). Assume 
\begin{equation}\label{eq:ass_casei}
 F(\bX\mymax{n})-F(\bXnmin) \leq \frac{\gam}{L}\big(1-\clam\big)\clam|\nabla F(\bX\mymax{n})|^2.
\end{equation}
Appealing to the descent property \eqref{eq:heaviest} for the heaviest agent 
$i=i^+_{n+1}$, then
\begin{equation}\label{eq:telesc}
\begin{split}
F(\bXnplusmin) & \leq F(\bx^{n+1}_{i^+_{n+1}}) \leq
F(\bX\mymax{n}) -\frac{2\gam}{L}\big(1-\clam\big)\clam|\nabla F(\bX\mymax{n})|^2 \\
 & \leq F(\bXnmin) -\frac{\gam}{L}\big(1-\clam\big)\clam|\nabla F(\bX\mymax{n})|^2.
\end{split}
\end{equation}
The inequality on the left follows since $\bXnplusmin$ is  the global minimizer at $t^{n+1}$; the middle inequality is the  descent property  for the heaviest agent, \eqref{eq:heaviest}, 
and  the last inequality follows from the assumed bound \eqref{eq:ass_casei}.\newline
Case (ii). We remain with the case
\begin{equation}\label{eq:ass_caseii}
F(\bX\mymax{n})-F(\bXnmin) \geq \frac{\gam}{L}\big(1-\clam\big)\clam|\nabla F(\bX\mymax{n})|^2.
\end{equation}
In this case, the (parent of the) heaviest agent positioned at $\bX\mymax{n}=\bx^n_{i^+_{n+1}}$, must different from the minimizer positioned at $\bXnmin=\bx^n_{i^-_n}$, or else $F(\bX\mymax{n})=F(\bXnmin)$ which is already covered by case (i),  \eqref{eq:ass_casei}.  Therefore, the heaviest agent had to shed a fraction of its mass, 
$\eta_+^n m^n_{i^+_{n+1}}$, which will be transferred to that minimzier:
\[
\left\{
\begin{split}
m^{n+1}_+ & = m^{n}_+ - \eta_+^n m^{n}_+, \qquad  \qquad m^{n}_+ := m^n_{i^+_{n+1}}, \ \ m^{n+1}_+ := m^{n+1}_{i^+_{n+1}}, \ \ \eta_+^n = \frac{F(\bX\mymax{n})-F(\bXnmin)}{F^{{}^n}_{\text{max}}-F^{{}^n}_{\text{min}}}, \\
 m^{n+1}_- & = m^n_{i^-_n} + \eta_+^n m^{n}_+ + (\textnormal{mass from other heavier agents}) \ldots, \quad \qquad m^{n+1}_- := m^{n+1}_{i^-_{n}}
 \end{split}\right.
 \]   
It follows that  the \emph{relative} mass of that minimizer is at least as large as
\[
m^{n+1}_- >\eta\mymax{n}m\mymax{n}= \frac{\eta\mymax{n}}{1-\eta\mymax{n}}m^{n+1}_+ \ \ \leadsto \ \ \relminin= \frac{m^{n+1}_-}{m^{n+1}_+} > \frac{\eta\mymax{n}}{1-\eta\mymax{n}}.
\] 
Recall that  the transition factors, $\etain$, are determined by the relative heights, \eqref{eq:etai}. Using the assumed bound of case (ii), \eqref{eq:ass_caseii}, we conclude
\[
\relminin >\eta_+^n = \frac{F(\bX\mymax{n})-F(\bXnmin)}{F^{{}^n}_{\text{max}}-F^{{}^n}_{\text{min}}}\geq \frac{\gam}{ML}\big(1-\clam\big)\clam|\nabla F(\bX\mymax{n})|^2.
\]
As before, the descent property \eqref{eq:allis} for $\bXnmin=\bx^n_{i^-_n}$  together with the lower bound we secured for  $\relminin$ in this case, imply
\begin{equation}\label{eq:telesb}
\begin{split}
F(\bXnplusmin)     \leq F(\bx^{n+1}_{i^-_n}) &\leq F(\bx^n_{i^-_n}) - \frac{2\gam}{L}\big(1-\clam\relminin\big)\clam\relminin|\nabla F(\bx^n_{i^-_n})|^2
 \\
& \leq F(\bXnmin)-\frac{2\gam^2}{ML^2}(1-\clam)^2\clam^2|\nabla F(\bX\mymax{n})|^2\times|\nabla F(\bXnmin)|^2.
\end{split}
\end{equation}
Combining   \eqref{eq:telesc} and \eqref{eq:telesb}, we find 
\begin{equation}\label{eq:telesf}
\begin{split}
F(\bXnplusmin) & \leq F(\bXnmin) -\frac{1}{C} \min\big\{|\nabla F(\bX\mymax{n})|^{2},   |\nabla F(\bX\mymax{n})|^{2}\times|\nabla F(\bXnmin)|^{2}\big\}\\
 & \leq F(\bXnmin) -\frac{1}{C}|\nabla F(\bX\mymax{n})|^{2} \times \min\big\{1,|\nabla F(\bXnmin)|^{2}\big\},
\end{split}
\end{equation}
with 
\begin{equation}\label{eq:telesC}
 C=\max\Big\{\frac{L}{\gamma(1-\clam)\clam},\frac{ML^2}{2\gamma^2(1-\clam)^2\clam^2}\Big\},
 \end{equation}
and a telescoping sum of \eqref{eq:telesC} implies the desired bound \myr{on the right of} \eqref{eq:summability} \myr{with $F(\bX_-^0)-F(\bx^*) \leq M$}. $\square$

\begin{remark}
We observe that the summability bound \eqref{eq:summability} is driven by a worst case scenario  alluded in case (ii) above. In this case, there is a potentially large difference of heights between the minimizing agent and heaviest agent; consequently, the descent property of the relatively lighter minimizer could be small and we had to rely on the descent property of the heaviest agent in \eqref{eq:telesc}, which led to \eqref{eq:telesb}.
\end{remark}
\smallskip
The descent rates of different agents can be arbitrarily slow,  due to  their time-dependent mass. The  summability bound \eqref{eq:summability} depends solely on the time sequence of  SBGD minimzers, $\{\bXnmin\}$, and (the parents of the) heaviest agents, $\{\bX\mymax{n}\}$, but it is  independent of the lightweight agents. 
Eventually, for large enough $n$, the minimizers and heaviest agents of SBGD${}_{pq}$ coincide into one time sequence, $\{\bX^n\}$.
The key point is that time sub-sequences, $\{\bX^{n_\alpha}\}$, satisfy a Palais-Smale condition \cite[{\S}II.2]{struwe2000variational}: by monotonicity, $F(\bX^{n_\alpha})\leq \max_i F(\bx^0_i)$ while $\nabla F(\bX^{n_\alpha}) \stackrel{\alpha\rightarrow \infty}{\longrightarrow}0$.

\begin{theorem}\label{thm:SBGD convergence}
Consider the loss function $F\in C^2(\Omega)$ such that the Lip bound \eqref{eq:FisC2} holds and let $\{\bXnmin\}_{n\geq0}$ denote the time sequence of SBGD  minimizers, \eqref{eq:SBGD},\eqref{eq:step length}.   Then $\{\bXnmin\}_{n\geq0}$ consists of one or more sub-sequences, $\{\bXnmina, \ \alpha=1,2,\ldots,\}$, that converge to a band of local minima with equal heights,
\[
\bXnmina \stackrel{n_\alpha\rightarrow \infty}{\longrightarrow}\bX^*_\alpha \ \ \textnormal{such that} \  \nabla F(\bX^*_\alpha)=0, \ \textnormal{and} \ F(\bX^*_\alpha)=F(\bX^*_\beta)
\]
In particular, if $F$ admits only distinct local minima in $\Omega$ (i.e., different local minima have different heights), then the whole sequence $\bX^n$ converges to a minimum.
\end{theorem}
\noindent \emph{Proof.} Since we assume  the sequence $\{\bXnmin\}$ is bounded in $\Omega$, it has  a converging sub-sequences. Take \emph{any} such converging sub-sequence $\bXnmina \rightarrow \bX^*_\alpha \in \Omega$. By \eqref{eq:summability}, $\nabla F(\bXnmina) \rightarrow 0$ for all sub-sequences, and hence $\bX^*_\alpha$ are local minimizers, $\nabla F(\bX^*_\alpha)=0$. Moreover, since $F(\bXnmin)$ is a decreasing, 
all $F(\bX^*_\alpha)$ must have the same `height'. The collection of equi-height minimizers  $\{\bX^*_\alpha \ \big| \  F(\bX^*_\alpha)=F(\bX^*_\beta)\}$ is the limit-set of $\{\bXnmin\}$. $\square$
\begin{remark} \myr{We  reiterate our earlier comment that the convergence analysis does not take into account the finitely many agents in the `survival of the fittest' policy. Putting efficiency aside, theorem \ref{thm:SBGD convergence} applies  to the general case where  all agents from initial crowd  survive the SBGD iterations.}
\end{remark}
\subsection{Flatness and convergence rate}
Theorem \ref{thm:SBGD convergence} indicates the convergence of SBGD without imposing any convexity condition on the loss function $F$, and therefore it comes without any rate. To quantify convergence \emph{rate}, we need  access to the fact that $F$ should be `curved up', at least  within a sufficiently small neighborhood of a local minimum $\bX^*_\alpha$. In the simplest case, $F$ may be assumed to be locally convex. However, one must take into account that $F$ may be more flat than just quadratic convexity. Indeed, these relatively flat local minima are the main hurdle in non-convex optimization. A precise classification for the level of `flatness' is offered by the Lojasiewicz condition.
Accortoing \myr{to} Lojasiewicz inequality, \cite{law1965ensembles,lojasiewicz1993geometrie}, if 
$F$ is analytic in $\Omega$ then  for every critical point of F, $\bxmin\in\Omega$, there exists a neighborhood $\calN_* \ni\bxmin$ surrounding $\bxmin$, an exponent $\beta\in (1,2]$ and a constant $\mu>0$ such that
\begin{equation}\label{eq:Lojasiewicz}
\mu|F(\bx)-F(\bxmin)| \leq |\nabla F(\bx)|^\beta, \qquad \forall \bx\in \calN_*.
\end{equation}
The exponent $\beta$ is tied to the \emph{flatness} of $\nabla F$ at $\bx=\bxmin$: if $\nabla F(\bx)$ vanishes of order $m$ at $\bx=\bxmin$, then $\displaystyle \beta=\frac{m+1}{m}$. In the particular case of local convexity,  $\bxmin$ is a simple minimum and  \eqref{eq:Lojasiewicz} is reduced to the \myr{Polyak}-Lojasiewicz condition \cite{polyak1964gradient} corresponding to $\beta=2$
\begin{equation}\label{eq:PL}
\mu\big(F(\bx)-F(\bxmin)\big) \leq |\nabla F(\bx)|^2, \qquad \forall \bx\in \calN(\bxmin).
\end{equation}
A smaller value of $\beta<2$ indicates a more flat configuration of  $F$ in a region $\calN_* \ni \bxmin$. \newline
In the theorem below, we restrict attention to SBGD${}_{1,1}$ and we assume that $n$ is large enough which allows us to treat only  the canonical scenario 
where minimizers and heaviest agents coincide, $\bX\mymax{n}=\bXnmin$.
\begin{theorem}\label{thm:GD-backtracking with PL}
 Consider the loss function $F\in C^2({\mathcal C})$ such that the Lip bound \eqref{eq:FisC2} holds, with minimal flatness $\beta$.  Let $\{\bXnmin\}_{n\geq0}$ denote the time sequence of SBGD minimizers,  \eqref{eq:SBGD},\eqref{eq:step length}.  Then, there exists a constant, $C=C(\gamma,\lambda,\mu)$, such that 
\begin{equation}\label{eq:GD PL convergence}
    F(\bXnmina)-F(\bX^*_\alpha) \left\{\begin{array}{ll}
   \displaystyle \leq \Big (1-\frac{2\mu \gamma\clam (1-\clam)}{L}\Big)^n\big(\min_i F(\bx^0_i)-F(\bxmin)\big), & \beta=2\\ \\
     \displaystyle \lesssim C \Big(\frac{1}{n_\alpha}\Big)^{\frac{\beta}{2-\beta}}, & \beta\in (1,2) \end{array}\right.
 \end{equation}
\end{theorem}
Observe that as `flatness', increases, $\beta$ is decreasing and the exponential decay in \eqref{eq:GD PL convergence}${}_1$ is replaced by a polynomial decay which may slow down all the way to first-order decay, $\nicefrac{1}{n_\alpha}$. 

\noindent 
\emph{Proof.} We limit ourselves to the canonical scenario in which the heaviest agent at $\bX\mymax{n}$ coincides with the minimizer at  $\bXnmin$. The descent property  for  the heaviest agent \eqref{eq:telesc} then reads
\[
F(\bXnplusmin)\leq F(\bXnmin) -\clam h_-|\nabla F(\bXnmin)|^2, \qquad  h_-:= \frac{2\gamma}{L}(1-\clam).
\]
We focus on the converging sub-sequence $\{\bXnmina\}$,
\[
F(\bXnplusmina)\leq F(\bXnmina) -\clam h_-|\nabla F(\bXnmina)|^2.
\]
Let us first discuss the quadratic case, $\beta=2$,  of  Polyak-Lojasiewicz condition \eqref{eq:PL}, which yields 
 \[
 F(\bXnplusmina) \leq F(\bXnmina)-\mu\lambda h_-\big(F(\bXnmina)-F(\bX^*_\alpha)\big), \qquad \bXnmina\in \calN_\alpha.
  \]
Rearranging we find
 \begin{equation}\label{eq:PL-gradient}
 F(\bXnplusmina)-F(\bX^*_\alpha)\leq \left(1-\mu\lambda h_-\right)\big(F(\bXnmina)-F(\bX^*_\alpha)\big),
 \end{equation}
which yields  exponential rate, \cite{polyak1964gradient,karimi2016linear}
 \[
F(\bXnmina)-F(\bX^*_\alpha)\leq \left(1-\mu \clam h_-\right)^n\big(F(\bx^0)-F(\bX^*_\alpha)\big).
\]
The case of general Lojasiewicz bound \eqref{eq:Lojasiewicz}  with $\beta<2$ yields
 that the error, $E_{n_\alpha}:=F(\bXnmina)-F(\bX^*_\alpha)$, satisfies 
\[
   E_{n_\alpha+1} \leq E_{n_\alpha}-\clam h_-(\mu E_{n_\alpha})^{\nicefrac{2}{\beta}},\qquad 
\bXnmina \in \calN_\alpha.
\]
The solution of this Riccati inequality  (e.g.,  \cite[Theorem 3.1]{tadmor1984large} for the limiting case $\beta=1$), yields
\[
F(\bXnmina)-F(\bX^*_\alpha)\lesssim_\mu \left\{|\min_i F(\bx^0_i)-F(\bX^*_\alpha)|^{-\nicefrac{1}{\beta'}}+ \clam h_-\mu^{\nicefrac{2}{\beta}} n_\alpha\right\}^{-\beta'}, \quad \beta'=\frac{\beta}{2-\beta}>1.
\]
and \eqref{eq:GD PL convergence} follows with $\ds C= \Big(\frac{L}{2 \gamma\clam (1-\clam)}\Big)^{\frac{\beta}{2-\beta}}(\nicefrac{1}{\mu})^{\frac{2}{2-\beta}}$. $\square$

\section{Numerical results --- one dimensional problems}\label{sec:results-1D}

\begin{figure}[htb]
    \centering
    \begin{subfigure}{0.4\textwidth}
    \centering
    \includegraphics[scale = 0.3]{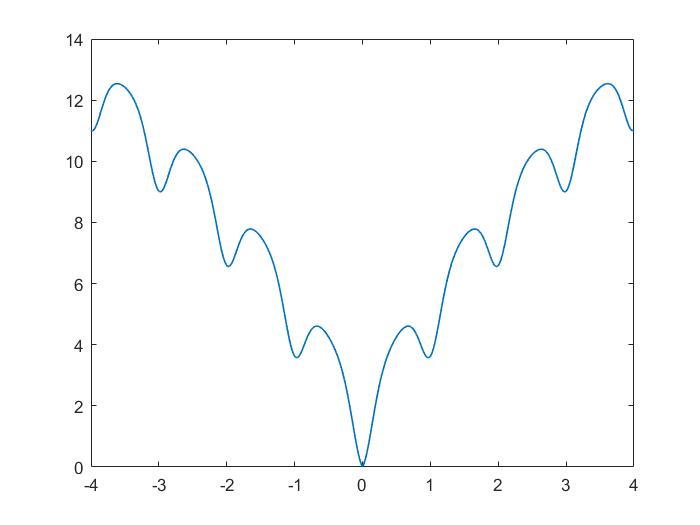}
    \subcaption{Ackley}
    \end{subfigure}
    \quad
    \begin{subfigure}{0.4\textwidth}
    \centering
    \includegraphics[scale = 0.3]{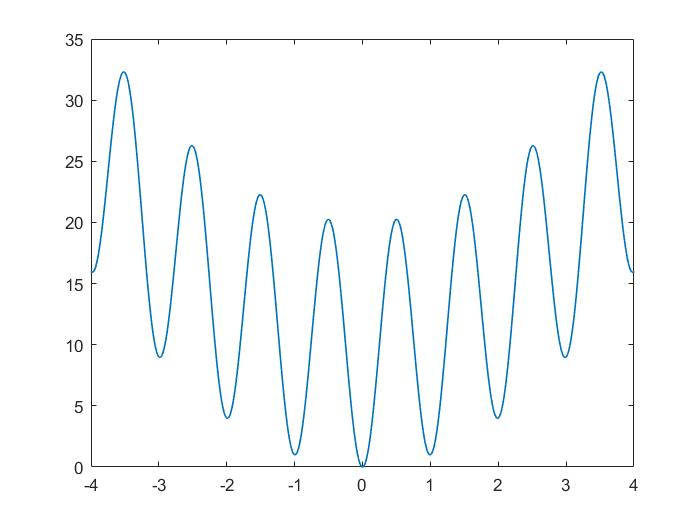}
    \subcaption{Rastrigin}
    \end{subfigure}
    \caption{Benchmark functions}
    \label{fig:benchmark functions}
\end{figure}

We use the swarm-based gradient descent method to search for the global minimizers of the 1D functions 

\[
\begin{split}
    & \textnormal{Ackley function}: \quad   \quad F_{\textnormal{Ackley}}(x) = -20 e^{-0.2|x_B|}- e^{\cos(2\pi (x_B))}+20+e+C, \ \ x_B:=x-B\\
    & \textnormal{Rastrigin  function}:  \quad F_{\textnormal{Rstrgn}}(x) =(x_B)^2-10\cos(2\pi(x_B))+10+C, 
    \end{split}
\]
with shift parameters $B,C$. Figure \ref{fig:benchmark functions} shows that both, the Ackley  and the Rastrigin function attain multiple local minimizers.

We implemented  the `vanilla' version of SBGD method $(p,q)=(1,1)$ (i.e., $\phi_p(\eta)=\eta$ and $\psi_q(\widetilde{m})=\widetilde{m}$), and with  parameters 
$\clam = 0.2, \  \gam = 0.9, \ h_0 = 1$. 
Here $\clam, \gam$ are the descent and shrinkage  parameters, and $h_0$ is the initial step length for backtracking line search. We record the results of   $m=200$ independent simulations of Algorithm \ref{alg:SBGD} with initial positions of the agents  uniformly distributed in the interval $[-3, 3]$.

\begin{figure}[h!]
    \centering
     \includegraphics[scale = 0.42]{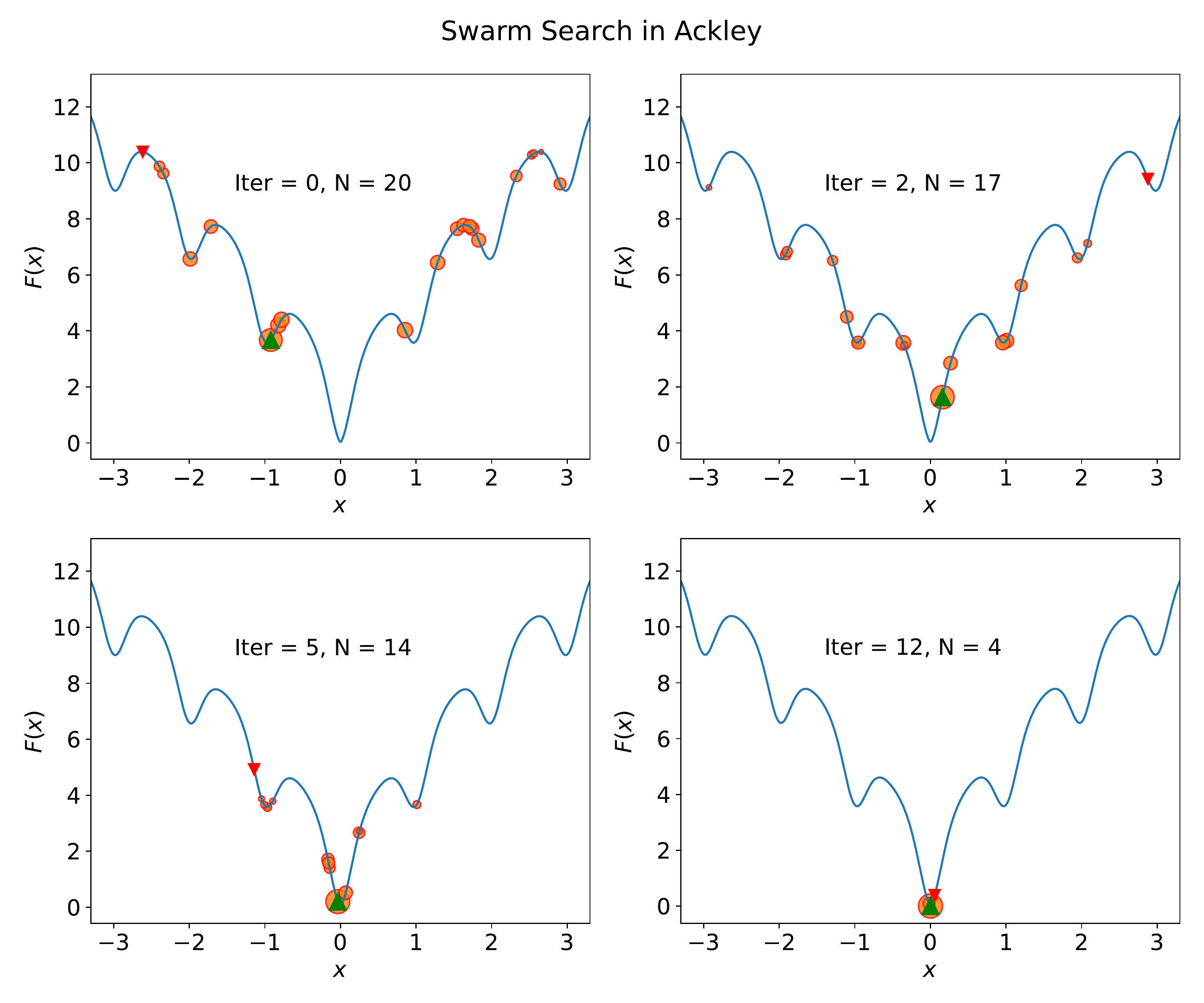}
     \caption{1D Ackley with $B = C = 0$.
    \myr{Four iterations of the SBGD visualized on the Ackley landscape show the dynamics of merged agents and convergence patterns.}}
    \label{fig:1d Ackley one exp.}
\end{figure}

\begin{figure}[h!]
       \centering
    \includegraphics[scale = 0.42]{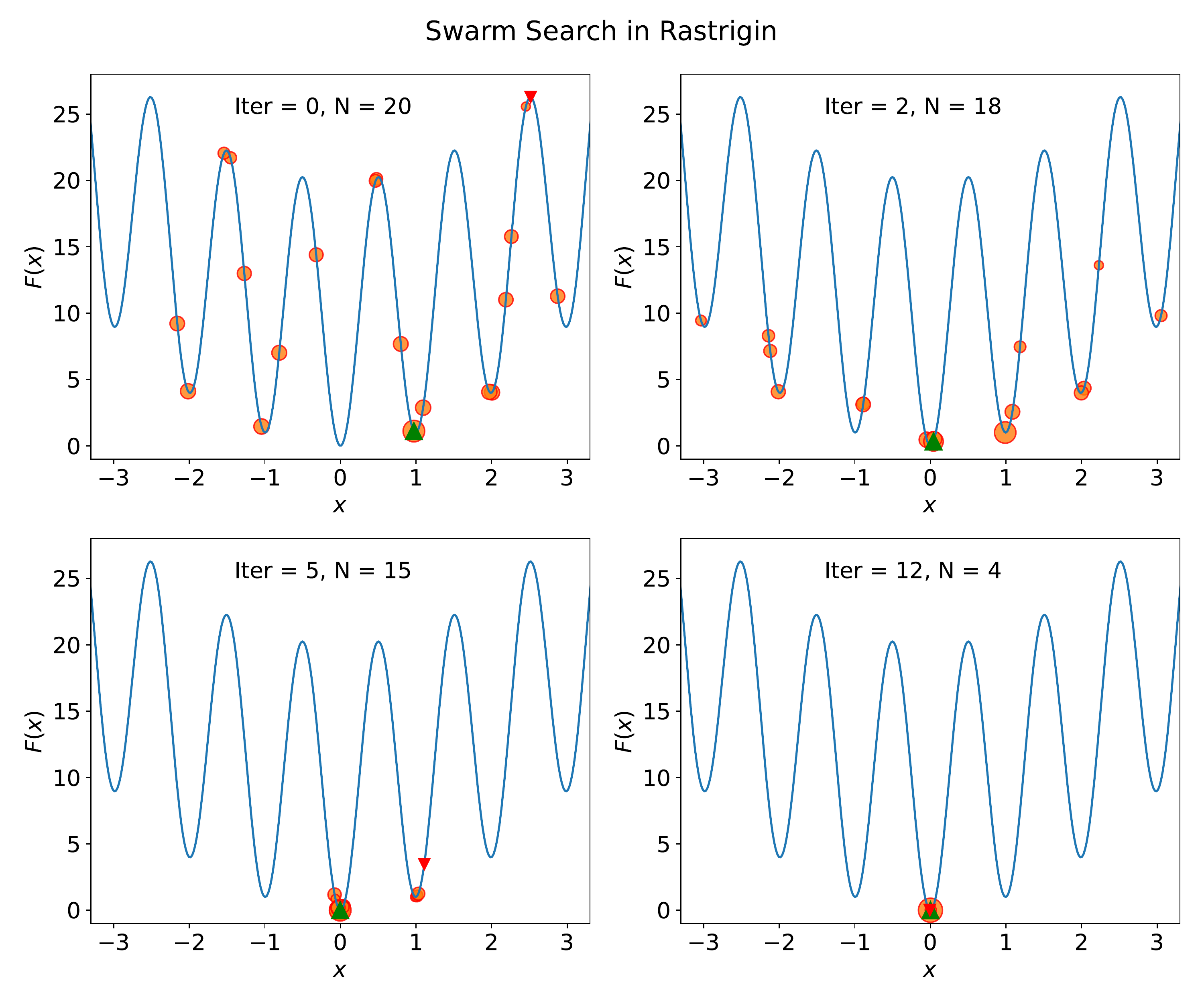}
     \caption{1D Rastrigin with $B = C = 0$.
    \myr{Four iterations of SBGD visualized on the Rastrigin landscape show show the dynamics of merged agents and convergence patterns.}}
    \label{fig:1d Rastrigin one exp.}
\end{figure}

To illustrate the behavior of the algorithm, we present the evolution of agents in one representative simulation for the Ackley function and the Rastrigin function in Figures \ref{fig:1d Ackley one exp.} and \ref{fig:1d Rastrigin one exp.}, where $N = 20$ agents are applied to search the domain.  The blue line depicts the target function, the red circles represent the agents. It is seen that as the iterations progress, the agents flock towards different minima. The agents getting stuck at local minimizers are gradually removed due to the mass transition, while those approaching the global minimizer are eventually merged into one agent. Table \ref{tab:1d Ackley Rastrigin} shows the results for the Ackley and the Rastrigin functions with a varying number of agents. The algorithm gives accurate solutions in all the test cases. With more agents applied the quality of the approximation is improved. However, we notice that improvement is significant as $N$ increases from 5 to 10, whereas increasing the number of agents from 10 to 20 does not result in much enhancement. 

\vspace{0.2in}
\begin{table}[htb]
\setstretch{1.5}
    \centering
    \begin{tabular}{|c|c||ccc|}
    \hline
                 B=0 &  & N=5 & N=10 & N=20\\ [0.5ex]
         \hline
    Ackley function $F_{\textnormal{Ackley}}$& \mbox{success rate} & $99.50\%$ &  $100\%$ & $100\%$\\
     & $\displaystyle \oneoverd \EE|x_{SOL}-x^*|^{2}$& $4.70e^{-3}$ & $8.21e^{-10}$ & $6.89e^{-10}$ \\ [0.5ex] 
     \hline
     Rastrigin function $F_{\textnormal{Rstrgn}}$& \mbox{success rate} & $97.00\%$ & $100\%$ & $100\%$\\
      & $\displaystyle \oneoverd \EE|x_{SOL}-x^*|^{2}$ & $4.45e^{-2}$ & $5.18e^{-10}$ & $1.89e^{-10}$\\ [0.5ex]
      \hline
    \end{tabular}
    \medskip
\caption{Results of SBGD for 1D Ackley and Rastrigin functions, $m = 200$.}\label{tab:1d Ackley Rastrigin}
\end{table}
\vspace{0.02in}

\begin{figure}[h!]
    \centering
    \begin{subfigure}{0.4\textwidth}
    \includegraphics[scale = 0.3]{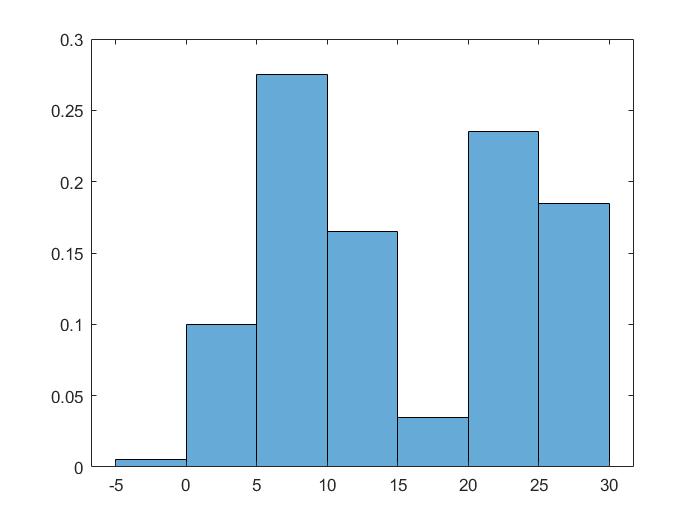}   
    \subcaption{SBGD, $N = 10$}
    \end{subfigure}
    \quad
    \begin{subfigure}{0.4\textwidth}
    \includegraphics[scale = 0.3]{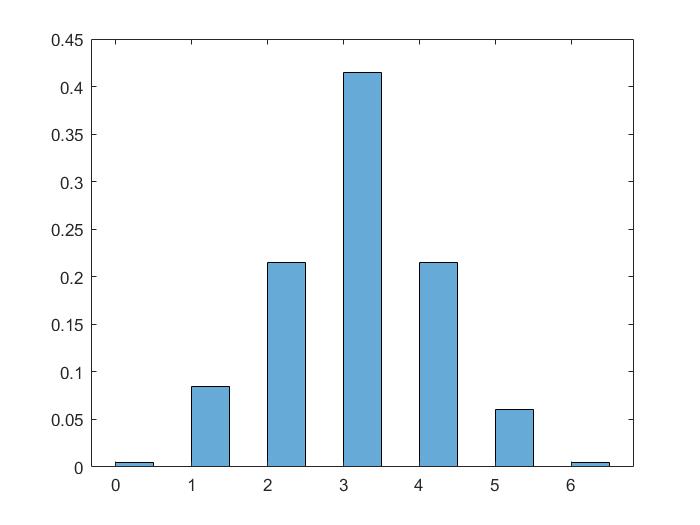}    
    \subcaption{\MAGD, $N = 10$}
    \end{subfigure}
    \\
    \begin{subfigure}{0.4\textwidth}
    \includegraphics[scale = 0.3]{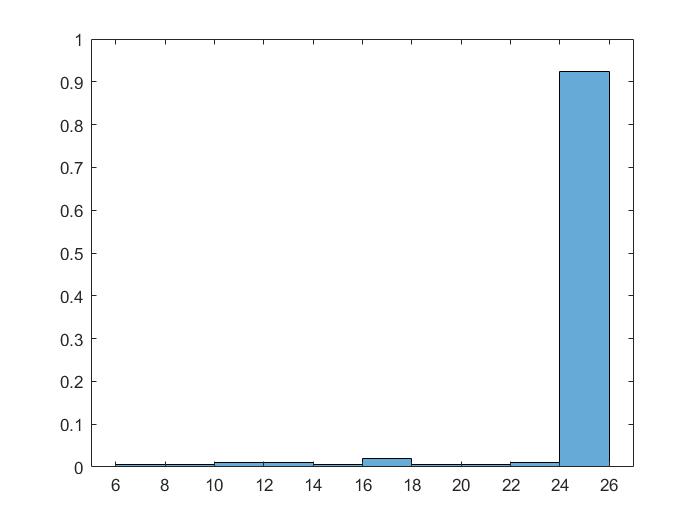}   
    \subcaption{SBGD, $N = 20$}
    \end{subfigure}
    \quad
    \begin{subfigure}{0.4\textwidth}
    \includegraphics[scale = 0.3]{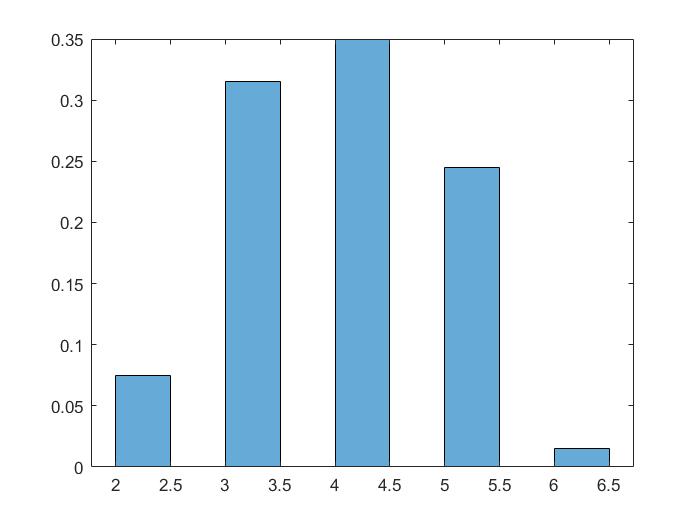}    
    \subcaption{\MAGD, $N = 20$}
    \end{subfigure}
    \\
    \begin{subfigure}{0.4\textwidth}
    \includegraphics[scale = 0.3]{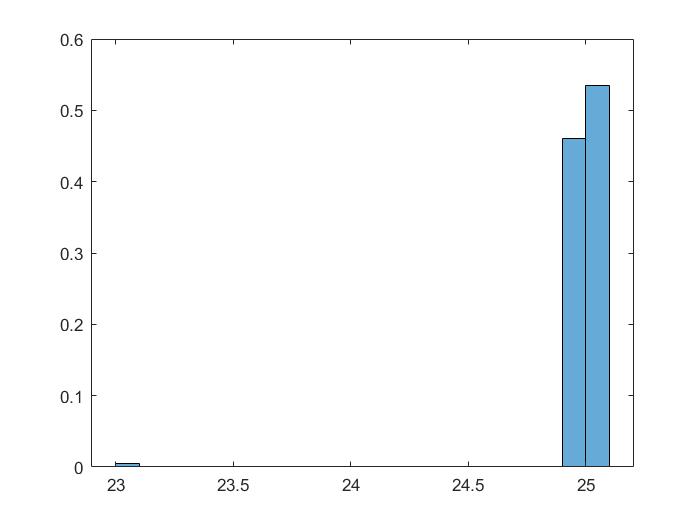}   
    \subcaption{SBGD, $N = 20$}
    \end{subfigure}
    \quad
    \begin{subfigure}{0.4\textwidth}
    \includegraphics[scale = 0.3]{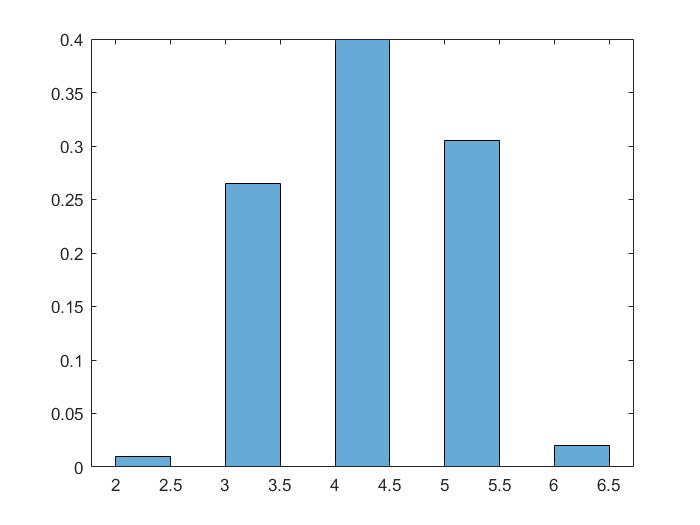}    
    \subcaption{\MAGD, $N = 20$}
    \end{subfigure}    
    \caption{Histograms of the shifted 1D Ackley function by 200 experiments. $B = 25$, $C = 5$. Global minimum $\min F_{\textnormal{Ackley}} = 5$ is attained at $\bx_{*} = 25$.}
    \label{fig:shifted 1d Ackley histogram}
\end{figure}

The key feature of the SBGD algorithm is communication,  reflected by the mass transitions between the agents. Such a mechanism dynamically adjusts the search strategies of different individuals. It is of interest to verify the benefit of  communication. We compare the results by the SBGD method with those by the non-communicating multi-agents gradient descent, GD(BT), of which all the agents conduct gradient descent search independently. Tables \ref{tab:shifted 1D Ackley SBGD} and \ref{tab:shifted 1D Ackley MGD}   present the results by SBGD and \MAGD methods obtained from $m = $200 experiments computed with varying shift parameter $B$ and different numbers of agents. We observe that when the initial distribution centers at $\bxmin$ or is moderately shifted, both SBGD and \MAGD methods are able to find the global minimizer. However, the advantage of SBGD  is more pronounced as $\bxmin$ is shifted farther away from the center of initial data, and the non-communicating \MAGD fails to give correct solutions since its exploration is restricted to the neighborhood of its initial data. Figure \ref{fig:shifted 1d Ackley histogram} displays the scenario  $B = 25$, where the initial distribution is strongly shifted away from the global minimizer. The \MAGD fails to find the global minimum. In contrast, the SBGD method employs light agents to conduct a more aggressive search in a larger area, and the algorithm ends at the correct minimizer at a surprisingly high success rate. 
\begin{remark}
Observe that when the initial data are  placed far from the global minmum,  the presence of a strong shift, $B=25$, requires sufficiently many swarming agents  $N> 20$, in order to secure  a  success rate $> 90\%$ and drive the  expected value of the error, $\EE|\bx_{SOL}-\bxmin|^{2} < 0.5$.
\end{remark}

\vspace{0.2in}
\begin{table}[htb]
    \begin{center}
    \begin{tabular}{|m{1.5cm}|m{3cm}|| m{1.8cm} m{1.8cm} m{1.8cm}|}
    \hline
          $\bxmin = B$  & & N=10 & N=20 & N=30\\ [0.5ex]
        \hline 
        $B = 0$  & \mbox{success rate} 
        $\displaystyle \oneoverd \EE|x_{SOL}-\xmin|^{2}$ & 
        $100\%$ $8.42e^{-10}$ & $100\%$ $8.37e^{-10}$ & 
        $100\%$ $3.38e^{-10}$\\  [0.5ex]
          \hline
        $B = 5$  & \mbox{success rate} 
        $\displaystyle \oneoverd \EE|x_{SOL}-\xmin|^{2}$ & 				
        $100\%$   $8.41e^{-10}$  & $100\%$ $7.58e^{-10}$ & 
        $100\%$ $5.01e^{-10}$ \\ [0.5ex]
            \hline 
        $B = 15$  & \mbox{success rate}
        $\displaystyle \oneoverd \EE|x_{SOL}-\xmin|^{2}$
         & $98.5\%$ $1.41e^{-2}$  & $100\%$ $4.69e^{-3}$ 
         & $100\%$ $8.27e^{-10}$\\ [0.5ex]
            \hline
        $B = 25$  & \mbox{success rate} 
        $\displaystyle \oneoverd \EE|x_{SOL}-\xmin|^{2}$
        & $45.5\%$ $1.48e^{+2}$ & $89.0\%$ $1.49e^{+1}$ & 
        $98.5\%$ $3.28e^{-1}$\\ [0.5ex]
        \hline 
    \end{tabular}
    \medskip
  \caption{Shifted 1D Ackley, results by SBGD, $m = 200$.}\label{tab:shifted 1D Ackley SBGD}
  \end{center}  
\end{table}
\vspace{0.2in}

\vspace{0.2in}
\begin{table}[htb]
    \begin{center}
    \begin{tabular}{|m{1.5cm}|m{3cm}|| m{1.8cm} m{1.8cm} m{1.8cm}|}
    \hline
          $\bx^* = B$  &  & N=10 & N=20 & N=30\\ [0.5ex]
        \hline
        $B = 0$  &  \mbox{success rate}  
        $\displaystyle \oneoverd \EE|x_{SOL}-x^*|^{2}$ &
        $100\%$ $8.60e^{-10}$ & $100\%$  $1.36e^{-9}$ & 
        $100\%$ $1.29e^{-9}$\\ [0.5ex]
       \hline
        $B = 5$  & \mbox{success rate} 
        $\displaystyle \oneoverd \EE|x_{SOL}-x^*|^{2}$ & 
        $100\%$ $8.51e^{-10}$ & $100\%$ $1.25e^{-9}$ & 
        $100\%$ $1.21e^{-9}$\\ [0.5ex]
           \hline
        $B = 15$  & \mbox{success rate} 
        $\displaystyle \oneoverd \EE|x_{SOL}-x^*|^{2}$ &
       $46.5\%$ $6.44e^{+1}$ & $75.0\%$ $2.26e^{+1}$ & 
        $85.5\%$ $1.14e^{+1}$\\ [0.5ex]
          \hline
        $B = 25$  & \mbox{success rate} 
        $\displaystyle \oneoverd \EE|x_{SOL}-x^*|^{2}$ & 
        $0\%$ $4.86e^{+2}$ & $0\%$ $4.50e^{+2}$  & 
        $0\%$ $4.38e^{+2}$\\ [0.5ex]
      \hline 
    \end{tabular}
    		\medskip
\caption{Shifted 1D Ackley, results by \MAGD (no communication), $m = 200$.}\label{tab:shifted 1D Ackley MGD}
\end{center}
\end{table}

\begin{figure}[h!]
    \centering
    \begin{subfigure}{0.4\textwidth}
    \includegraphics[scale = 0.3]{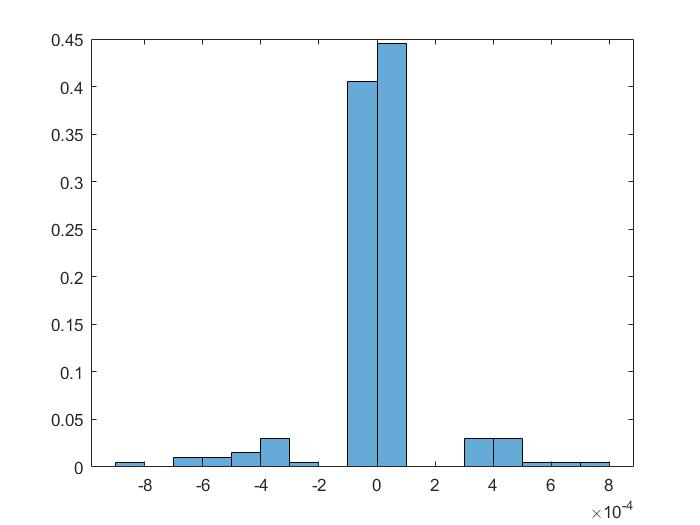}
    \subcaption{Ackley}
    \end{subfigure}
    \quad
    \begin{subfigure}{0.4\textwidth}
    \includegraphics[scale = 0.3]{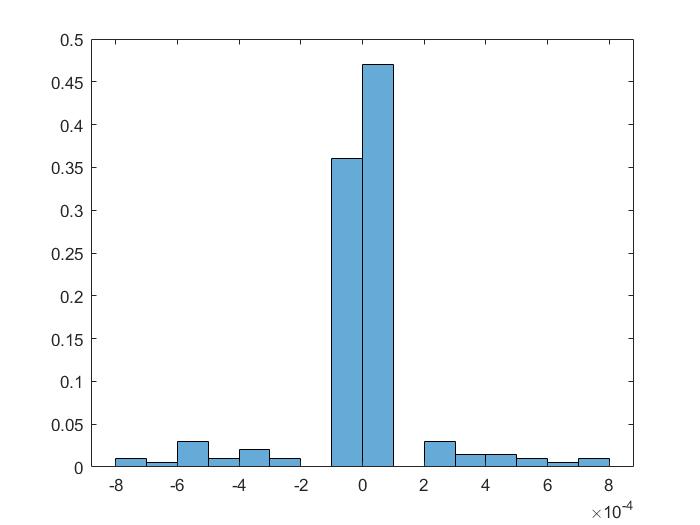}
    \subcaption{Rastrigin}
    \end{subfigure}
    \caption{Histograms of 1-D optimization, $B=C=0$,  200 simulations, bin width = $10^{-4}$.}
    \label{fig:1d histogram}
\end{figure}

\begin{figure}[h!]
    \centering
    \begin{subfigure}{0.45\textwidth}
    \centering
    \includegraphics[scale = 0.3]{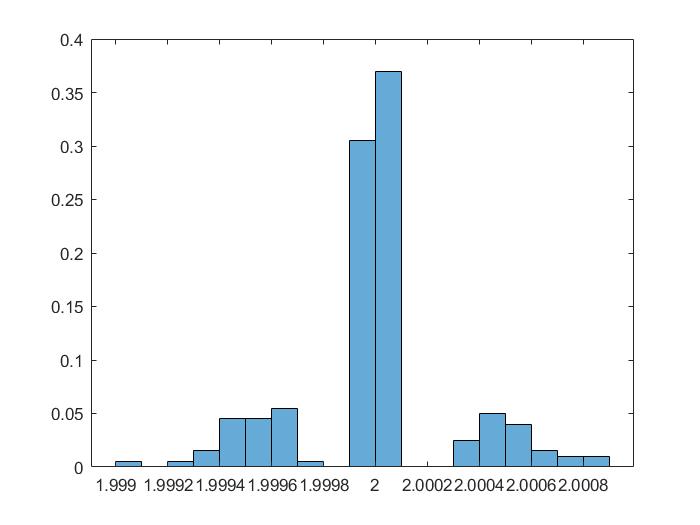}
    \subcaption{Ackley: $B = 2$, $C = 5$, $\argmin F_{A} = 2$.}
    \end{subfigure}
    \quad 
        \begin{subfigure}{0.45\textwidth}
    \centering
    \includegraphics[scale = 0.3]{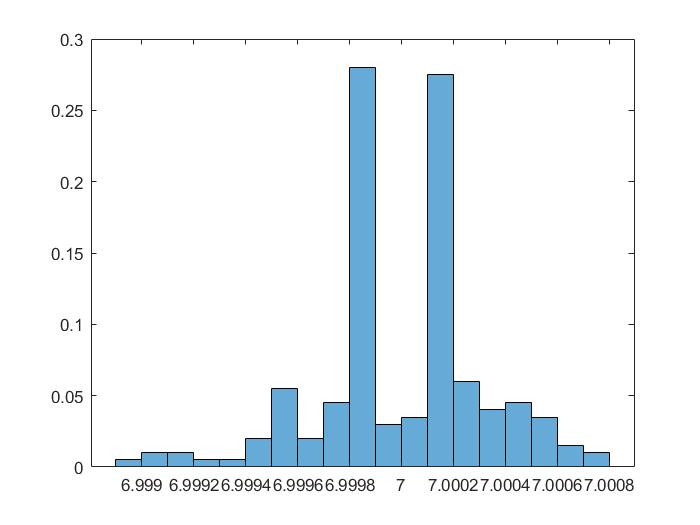}
    \subcaption{Ackley: $B = 7$, $C = 5$, $\argmin F_{A} = 7$.}
    \end{subfigure}
    \\
    \begin{subfigure}{0.45\textwidth}
    \centering
    \includegraphics[scale = 0.3]{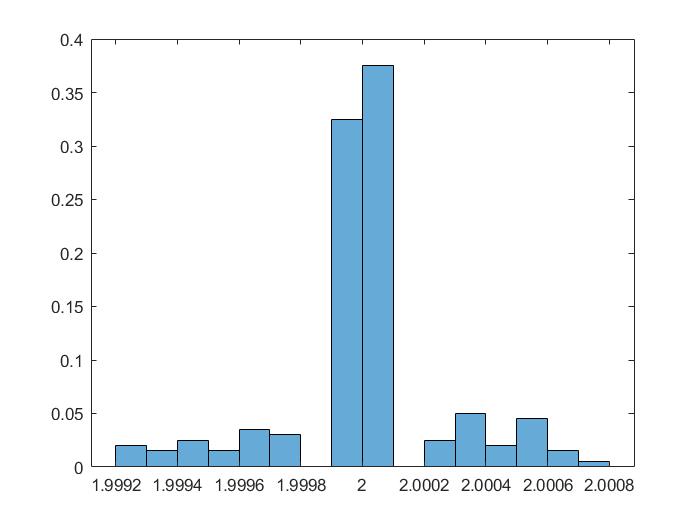} 
    \subcaption{Rastrigin: $B = 2$, $C = 5$, $\argmin F_{B} = 2$.}
    \end{subfigure}
    \quad 
        \begin{subfigure}{0.45\textwidth}
    \centering
    \includegraphics[scale = 0.3]{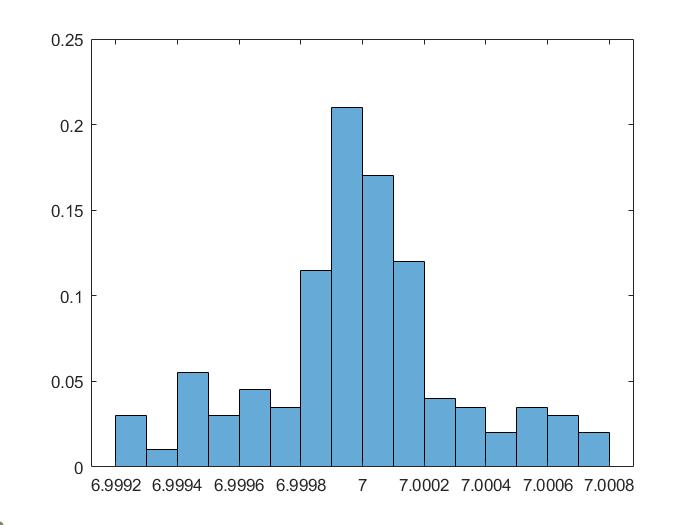}
    \subcaption{Rastrigin: $B = 5$, $C = 5$, $\argmin F_{B} = 7$.}
    \end{subfigure}
    \caption{Histograms of 1D optimization with shifting, 200 simulations, bin width = $10^{-4}$.}\label{fig:1d histogram shifted}
\end{figure}

\section{Numerical results --- two-dimensional problems}\label{sec:results-2D}

\subsection{Optimization of parameters} Let us first comment on the general issue of optimization of parameter space.
As noted in \cite{wilson2017marginal},  once the parameterization of an adaptive GD method is fixed,  it may not yield as good as or better results as simpler GD methods,  yet  adaptivity should be `judged' after being  optimized in space parameter, \cite{choi2019empirical}.  In this context, one can argue that optimizing a single agent method in parameter space is equivalent to a selective choice among many simulations of non-communicating multi-agent dynamics, whereas the swarm-based approach provides a \emph{dynamic}, `on the fly' selection of optimized parameters, which is precisely the type of comparisons we make below. 

On the other hand, our SBGD method depends on several parameters: the initial step $h_0$, the descent parameter $\clam$ and shrinkage parameter $\gamma$ tied to the backtracking, the tolerance parameters \eqref{eq:tols} and the $(p,q)$ parameters. In the multi-dimensional computations reported below, we do \underline{not} optimize these SBGD parameters. 
Thus, unless otherwise stated, we examine the performance of the SBGD method  in Algorithm \ref{alg:SBGD} 
with initial step size $h_0=1$, a descent parameter $\clam=0.2$, a shrinkage parameter $\gamma=0.9$ and the threshold parameters \eqref{eq:tols}. We begin here with the `vanilla' version of SBGD, $(p,q)=(1,1)$, although later  we shall find out that  the choice $(p,q)=(2,1)$  seems  universally better. 
 We run  $m$  number of independent simulations, initiated with uniformly distributed positions, $\{\bx^0_i\}$ and  measure the success of the SBGD according to the proportion of its successful end results.

We illustrate the performance of the SBGD algorithm in multiple dimensions on  three benchmark test cases, \cite{benchmarks}. First, the Ackley function
\begin{equation}\label{eq:Ackley}
    F_{\textnormal{Ackley}}(\bx) = -20\exp\Big\{-\frac{0.2}{\sqrt{d}}\Big\{\sum^{d}_{i=1}\xBi^{2}\Big\}^{\nicefrac{1}{2}}\Big\}-\exp\Big\{\frac{1}{d}\sum^{d}_{i=1}\cos(2\pi \xBi)\Big\}+20+e+C.
\end{equation}
Second, the Rastrigin function
\begin{equation}\label{eq:Rastrigin}
    F_{\textnormal{Rstgin}}(\bx) = \frac{1}{d}\sum^{d}_{i=1}\Big\{\xBi^{2}-10\cos(2\pi \xBi)+10\Big\}+C.
\end{equation}
Here $d$ is the dimension of the ambient space, $\xB:=(x_1-B, \ldots, x_d-B)$ is the shifted variable in ${\mathbb R}^d$ and  $B, C\in\R$ are the shift parameters. Both functions attain the global minimum $C$ at the unique global minimizer $\bxmin = B$.

Third, we consider the   drop-wave function 
\begin{equation}\label{eq:drop-wave}
F_{\textnormal{Drop}}(\bx)= -\frac{1+\cos(12|\bx|)}{0.5|\bx|^2+2},
\end{equation}
with SBGD parameters $\clam = 0.3, \ \gam = 0.9, \ h_0 = 1$.
 The agents  are initialized with uniform distribution of positions in the hypercube $[-3, 3]^d$.

To evaluate the quality of the solution, we  make use of the \emph{success rate} among $m$ independent simulations. We consider a simulation to be successful if $ \bx_{SOL}$ is within the $d$-dimensional cube $ [\bxmin-0.25, \bxmin+0.25]^{d}$, e.g., \cite[\S4.2]{carrillo2021consensus}. This condition ensures that the approximate solution lies in the basin of attraction of the global minimizer. In fact, in a successful experiment the solution will lie in a much smaller neighborhood of $\bxmin$.

\subsection{SBGD compared with non-communicating GD(BT)}
\begin{figure}[h!]
    \centering
    \begin{subfigure}{0.4\textwidth}
    \includegraphics[scale = 0.3]{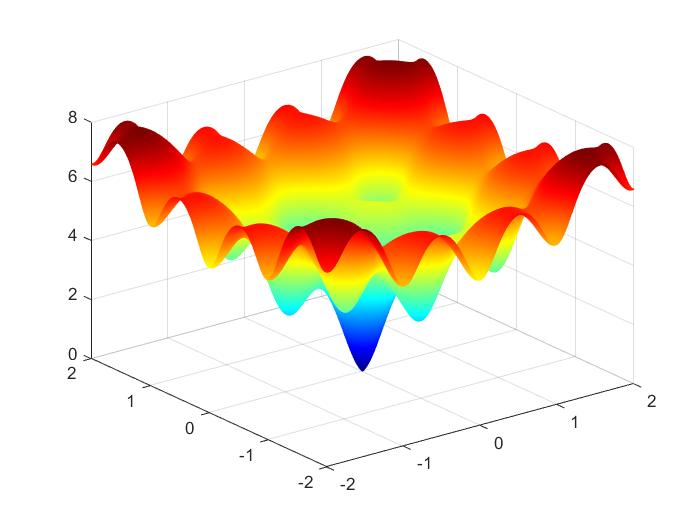}  
    \subcaption{2-D Ackley}
    \end{subfigure}
    \quad
    \begin{subfigure}{0.4\textwidth}
    \includegraphics[scale = 0.3]{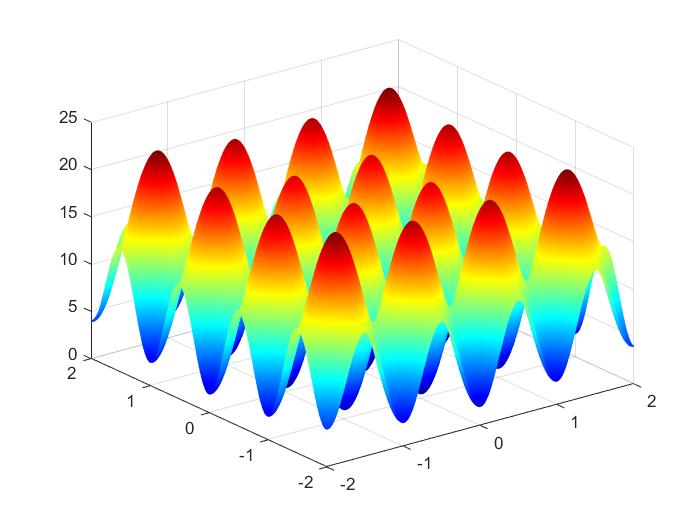}
    \subcaption{2-D Rastrigin }
    \end{subfigure}
    \quad
     \begin{subfigure}{0.4\textwidth}
    \includegraphics[scale = 0.3]{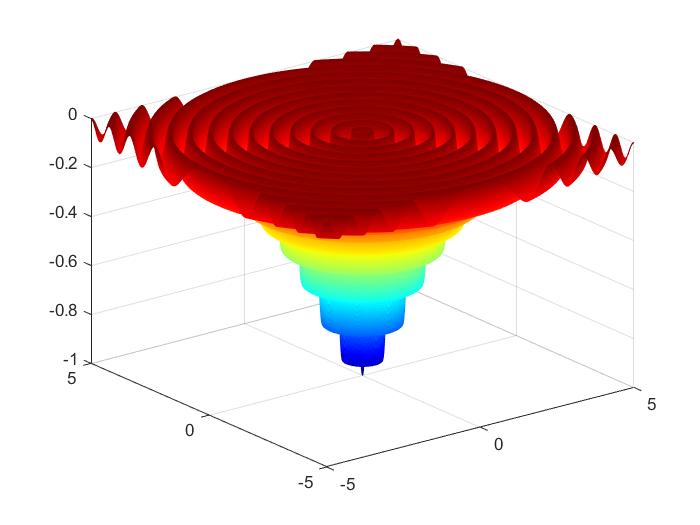}  
    \subcaption{2-D drop-wave function}
    \end{subfigure}
    \quad
    \begin{subfigure}{0.4\textwidth}
    \includegraphics[scale = 0.3]{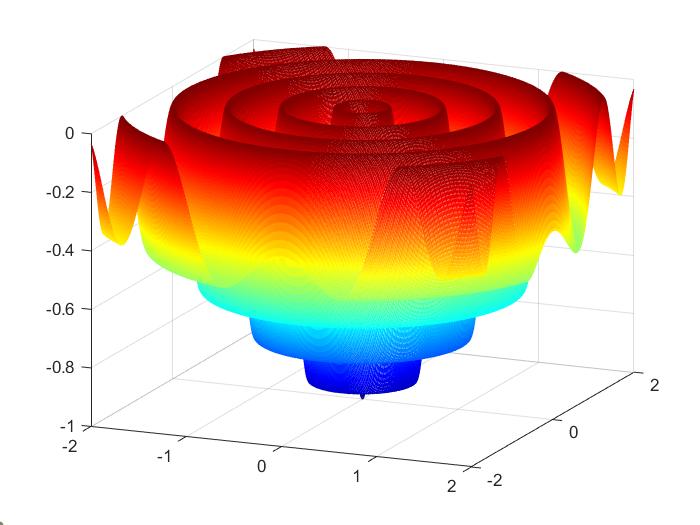}
    \subcaption{zoom into 2-D drop-wave function}
    \end{subfigure}
    \caption{Landscapes of benchmark functions in two dimensions}
    \label{fig:2D functions}
\end{figure}

We verify the advantage of the SBGD method in comparison to  the non-communicating  \MAGD algorithm. We consider the three benchmarks  of Ackley, Rastrigin and drop-wave functions  in  two dimensions. The landscapes are as shown in Figure \ref{fig:2D functions}. All functions have multiple local minimums. While the global minimum of the Ackley function is obviously lower than the other local minimums, the global minimum of the Rastrigin function is much less distinguishable. The drop-wave function has complex geometry with high frequency local minima and  sharp basins of attractions. Its global minimum $F(\bxmin) = -1$ is attained at the unique global minimizer $\bxmin = 0$. 

Table \ref{tab:2D Ackley SBGD vs GD(BT)} compares the success rates of SBGD and \MAGD methods in the Ackley test cases with varying function shifts $B$ and different numbers of agents $N$. It is observed that the two methods perform comparably well when the initial distribution is centered at $\bxmin$ or moderately shifted. But as the shift parameter $B$ increases, the SBGD method still achieves high success rates whereas the non-communicating \MAGD method fails. The results for the Rastrigin function, given in Table \ref{tab:2D Rastrigin SBGD vs  GD(BT)}, echo the same advantage of SBGD in shifted scenarios.

\begin{table}[h!]
\setstretch{1.5}
    \begin{center}
    \begin{tabular}{|m{1.5cm}|m{1.6cm}|| m{1.6cm} m{1.6cm} m{1.8cm}|}
    \hline
         $\bx^* = B$  & & N=25 & N=50 & N=100\\ [0.5ex]
        \hline
        $B = 0$  & SBGD \MAGD & 
        $98.0\%$ $100\%$ & $100\%$ $100\%$  & $100\%$ $100\%$ \\ [0.5ex]
             \hline
        $B = 5$  & SBGD \MAGD & 
        $93.6\%$ $71.2\%$ & $98.6\%$ $87.0\%$ & $99.8\%$ $99.2\%$ \\ [0.5ex]
           \hline
        $B = 10$  & SBGD \MAGD  & 
        $66.2\%$ $0\%$ & $90.8\%$ $0\%$ & $98.4\%$ $0.6\%$\\ [0.5ex]
            \hline 
         \end{tabular}
                \end{center}
                \medskip
       \caption{2D Ackley, success rates of SBGD vs. GD(BT), $m = 500$.}\label{tab:2D Ackley SBGD vs GD(BT)}
\end{table}

\begin{table}[h!]
\setstretch{1.5}
\vspace{0.2in}
    \begin{center}
    \begin{tabular}{|m{1.5cm}|m{1.6cm}|| m{1.6cm} m{1.6cm} m{1.8cm}|}
    \hline
      $\bx^* = B$  & & N=25 & N=50 & N=100\\ [0.5ex]
        \hline
        $B = 0$  & SBGD \MAGD & 
        $73.6\%$ $53.0\%$ & $95.2\%$ $76.4\%$ & $100\%$ $96.40\%$ \\ [0.5ex]
             \hline
        $B = 5$  & SBGD \MAGD  & 
        $44.4\%$ $15.8\%$  & $80.4\%$ $28.4\%$ & $99.2\%$ $56.4\%$\\ [0.5ex]
             \hline
        $B = 10$  & SBGD \MAGD & 
        $14.8\%$ $13.0\%$ & $44.0\%$ $30.4\%$ & $78.4\%$ $56.2\%$\\ [0.5ex]
          \hline 
    \end{tabular}
    \medskip
    \caption{2D Rastrigin, success rates of SBGD vs. GD(BT), $m = 500$.} \label{tab:2D Rastrigin SBGD vs GD(BT)}
    \end{center}
\end{table}

\begin{table}[h!]
\setstretch{1.5}
\vspace{0.2in}
    \begin{center}
    \begin{tabular}{|m{1.5cm}|| m{1.0cm} | m{1.0cm} | m{1.0cm}|}
    \hline
        & N=10 & N=20 & N=30\\ [0.5ex]
        \hline
          SBGD &         $90.5\%$ & $99.5\%$ & $100\%$  \\ [0.5ex]
             \hline
         \MAGD  & $15.0\%$ & $21.5\%$ & $35.5\%$\\ [0.5ex]
               \hline 
    \end{tabular}
    \medskip
    \caption{2D drop-wave, success rates of SBGD vs. GD(BT), $m = 500$.} \label{tab:2D drop-wave}
    \end{center}
\end{table}

\begin{table}[h!]
\setstretch{1.5}
\vspace{0.2in}
    \begin{center}
    \begin{tabular}{|m{2.5cm}|| m{1.5cm} | m{1.5cm} | m{1.5cm}|}
    \hline
        $\EE|\bx_{SOL}-\bxmin|^2 $ & $\clam=0.3$ & $\clam=0.5$ & 
        $\clam=0.7$\\ [0.5ex]
        \hline
          \ \ \ SBGD &     2.07e-1    & 6.18e-2 & 1.29e-1  \\ [0.5ex]
             \hline
         \ \ \ \MAGD  & 4.25e-2 & 6.08e-2 & 1.29e-1\\ [0.5ex]
               \hline 
    \end{tabular}
    \medskip
    \caption{Average error for 2D Rosenbrock function SBGD vs. GD(BT), $m = 500$ simulations equi-distributed at $[-4,-2]^2$.} \label{tab:2D SBGD for Rosenbrock}
    \end{center}
\end{table}

Table \ref{tab:2D drop-wave} compares the success rates of SBGD and \MAGD methods in the  test case of the 2D drop-wave function, \eqref{eq:drop-wave}.
We compare the success rates of SBGD with $\clam=0.3$ vs. the non-communicating \MAGD using $m=200$ simulations randomly initiated with uniform distribution at  $[-3, 3]^2$. The SBGD results reported in  Table \ref{tab:2D drop-wave} show a remarkable improvement over the performance of GD(BT); communication helps.

The SBGD does not always offer such a decisive lead over GD(BT). In table \ref{tab:2D SBGD for Rosenbrock} we compare the error of SBGD vs. \MAGD for the 2D Rosenbrock function
\begin{equation}\label{eq:Rosenbrock}
F_{\textnormal{Rsnbrk}}(\bx) = (1 - x_1)^2 + 100(x_2 -x_1^2)^2.
\end{equation}
The stiffness of $F_{\textnormal{Rsnbrk}}$ near its global minimum at $\bxmin = (1,1)$ produces  comparable results of the swarm dynamics  and non-communicating GD(BT). 

\subsection{SBGD${}_{pq}$ method --- dependence of $(p,q)$ parameters}\label{sec:pq}
We now turn to discuss the effect of adjusting the mass transition  with $\etainp$ and backtracking with $\relminq$.\newline
We implemented the SBGD${}_{pq}$ for shifted 2D Rastrigin
function \eqref{eq:Rastrigin} with shift $B=5$, and  drop-wave function \eqref{eq:drop-wave} using SBGD${}_{pq}$ with descent parameter $\clam=0.3$, using  $N$ agents  uniformly initialized at $ [-3, 3]^2$. The success rate $m=500$ simulations with different $(p,q)$ are recorded in Table \ref{tab:2D SBGDpq}.\newline
In both cases,  the results for  the `vanilla' SBGD, $(p,q)=(1,1)$, with a success rate of 80\% and respectively 90\%, are substantially improved to 97\%, once we use SBGD${}_{pq}$ with $(p,q)=(2,\nicefrac{1}{2})$. In this case, the use of $p=2$ enforces  more moderate mass transitions which seems to play a key role whenever  $F$ has steep basins of attractions, while $q=\nicefrac{1}{2}$ enforces a slower marching protocol for intermediate agents.

\begin{table}[h!]
\setstretch{1.5}
\vspace{0.1in}
    \begin{center}
    \parbox{.45\linewidth}
    {\begin{tabular}{|m{2.5cm} || m{1.2cm} | m{1.2cm} | m{1.2cm}|}
    \hline
       \backslashbox{$\phi_p(\eta)$}{$\psi_q(\widetilde{m})$}  & 
       \ \ $\sqrt{\widetilde{m}}$ & \ \ $\widetilde{m}$ & \ \ $\widetilde{m}^2$\\ [0.5ex]
        \hline
         \qquad $\sqrt{\eta}$ & 77.2\% & 43.0\% & 40.0\%\\ [0.5ex]
        \hline
		\qquad $\eta$    &   94.8\% & 80.4\% & 73.4\%\\ [0.5ex]
        \hline
		\qquad $\eta^2$ &  97.8\% & 94.0\% & 88.4\%\\ [0.5ex]
        \hline
		\qquad $\eta^{20}$ & 83.4\% & 95.4\% & 98.2\%\\ [0.5ex]
        \hline
    \end{tabular}
    \medskip
    \subcaption{2D Rastrigin function; $N=50$ agents} 
     }
    \medskip
      \qquad
        \parbox{.45\linewidth}
    {\begin{tabular}{|m{2.5cm} || m{1.2cm} | m{1.2cm} | m{1.2cm}|}
    \hline
       \backslashbox{$\phi_p(\eta)$}{$\psi_q(\widetilde{m})$}  & 
       \ \ $\sqrt{\widetilde{m}}$ & \ \ $\widetilde{m}$ & \ \ $\widetilde{m}^2$\\ [0.5ex]
        \hline
         \qquad $\sqrt{\eta}$ & 93.6\% & 83.0\%  & 79.8\%\\ [0.5ex]
        \hline
			 \qquad $\eta$    & 97.2\% & 90.4\% & 85.2\%\\ [0.5ex]
        \hline
		 \qquad $\eta^2$ & 97.6\% & 96.0\% & 90.2\% \\ [0.5ex]
        \hline
		 \qquad $\eta^{20}$ & 80.8\% & 98.8\% & 97.4\%\\ [0.5ex]
        \hline
     \end{tabular}
     \medskip
     \subcaption{2D Drop-wave function; $N=10$ agents}}
    \medskip
    \caption{SBGD${}_{pq}$ for shifted 2D Rastrigin and drop-wave with different $(p,q)$} 
    \label{tab:2D SBGDpq}
    \end{center}
\end{table}

\begin{table}[h!]
\setstretch{1.5}
\vspace{0.1in}
    \begin{center}
    \parbox{.47\linewidth}
    {\begin{tabular}{|m{2.5cm} ||  m{1.3cm} | m{1.3cm} | m{1.3cm}|}
    \hline
         \backslashbox{$\phi_p(\eta)$}{$\psi_q(\widetilde{m})$} & 
         \ \ $\sqrt{\widetilde{m}}$ & \ \ $\widetilde{m}$ & \ \ $\widetilde{m}^2$\\ [0.5ex]
        \hline
         \qquad $\eta$ & 2.53e-1 & 2.07e-1 & 1.78e-1\\ [0.5ex]
        \hline
			 \qquad $\eta^2$    & 8.51e-2 & 3.07e-1 & 1.70e-1\\ [0.5ex]
        \hline
		 \qquad $\eta^{20}$ & 6.39e-2 & 2.01e-1 & 2.74e-1 \\ [0.5ex]
        \hline
		 \qquad $\eta^{50}$ & 5.75e-2 & 1.88e-1 & 2.61e-1\\ [0.5ex]
        \hline
     \end{tabular}
     \medskip
     \subcaption{SBGD${}_{pq}$ with $\clam=0.3$.} 
     }
    \medskip
      \qquad
        \parbox{.47\linewidth}
    {\begin{tabular}{|m{2.5cm} ||  m{1.3cm} | m{1.3cm} | m{1.3cm}|}
    \hline
         \backslashbox{$\phi_p(\eta)$}{$\psi_q(\widetilde{m})$} & 
         \ \ $\sqrt{\widetilde{m}}$ & \ \ $\widetilde{m}$ & \ \ $\widetilde{m}^2$\\ [0.5ex]
        \hline
         \qquad $\eta$ & 6.79e-2 & 6.18e-2 & 6.05e-2\\ [0.5ex]
        \hline
			 \qquad $\eta^2$ &   4.97e-2 & 7.10e-2 & 6.08e-2\\ [0.5ex]
        \hline
		 \qquad $\eta^{20}$ & 4.65e-2 & 5.95e-2 & 6.74e-2 \\ [0.5ex]
        \hline
		 \qquad $\eta^{50}$ & 4.98e-2 & 6.02e-2 & 6.87e-2\\ [0.5ex]
        \hline
     \end{tabular}
     \medskip
     \subcaption{SBGD${}_{pq}$ with $\clam=0.5$.} 
     }
    \medskip
     \parbox{.47\linewidth}
    {\begin{tabular}{|m{2.5cm} ||  m{1.3cm} | m{1.3cm} | m{1.3cm}|}
    \hline
         \backslashbox{$\phi_p(\eta)$}{$\psi_q(\widetilde{m})$} & 
         \ \ $\sqrt{\widetilde{m}}$ & \ \ $\widetilde{m}$ & \ \ $\widetilde{m}^2$\\ [0.5ex]
        \hline
         \qquad $\eta$ & 7.11e-2 & 1.51e-1 & 1.29e-1\\ [0.5ex]
        \hline
			 \qquad $\eta^2$ &   6.83e-2 & 1.22e-1 & 1.29e-1\\ [0.5ex]
        \hline
		 \qquad $\eta^{20}$ & 1.21e-1 & 6.45e-2 & 1.48e-1\\ [0.5ex]
        \hline
		 \qquad $\eta^{50}$ & 1.26e-1 & 6.76e-2 & 1.37e-1\\ [0.5ex]
        \hline
     \end{tabular}
     \medskip
     \subcaption{SBGD${}_{pq}$ with $\clam=0.7$.} 
     }
    \medskip
    \caption{SBGD${}_{pq}$ simulation  of 2D Rosenbrock $F_{\textnormal{Rsnbrk}}$ with different $(p,q)$. A larger $\clam$ enforces a stronger descent property. The parameter $p=2$ is a most effective parameterization of SBGD${}_{pq}$.} 
    \label{tab:2D Rosenbrockpq}
    \end{center}
\end{table}

As a further  example we consider the 2D Rosenbrock function $F_{\textnormal{Rsnbrk}}$ in \eqref{eq:Rosenbrock}.
The difficulty arises when the iterations  approach the neighborhood of a global minimum $F(\bxmin) = 0$ attained at $\bxmin = (1,1)$. This  is due to the severe ``skew-ness'' of $F$  which is sensitive to the descent parameter $\clam$.  We employ $N = 30$ agents initialized with randomly distributed positions at $[-4,-2]^2$.   Success rates are computed among $m=500$ experiments. Table \ref{tab:2D Rosenbrockpq} records the \MAGD and SBGD${}_{pq}$ errors, measured by $\EE|\bx_{SOL} - \bxmin|^2$, for different values of descent parameter, $\clam$, and for different protocols of $(p,q)$. The sensitive dependence on $\clam$  is observed with both the non-communicating \MAGD and the SBGD${}_{pq}$.
Once again, the parameters $(p,q)=(2,\nicefrac{1}{2})$ 
yield the most stable performance of SBGD${}_{pq}$, whereas  other scaling
of \MAGD and SBGD are more sensitive to the choice of $\clam$.

In summary, Tables \ref{tab:2D SBGDpq} and \ref{tab:2D Rosenbrockpq}  indicate that while the results of SBGD${}_{pq}$ are mostly comparable, SBGD${}_{(2,\nicefrac{1}{2})}$ seems to provide optimal results, with main emphasize on $p=2$. 
At the same time, we conclude that the  tuning parameters, $(pq)$,  have a limited effect on the overall performance of SBGD${}_{pq}$ method. Accordingly, we did not optimize these parameters.\newline 
Motivated by these findings, we  focus below on two versions of SBGD (see also Tables \ref{tab:1D FB SBGD vs GD-BT} and \ref{tab:shifted 1D FB SBGD vs GD-BT}): the vanilla version,  SBGD${}_{1,1}$ with $\phi_p(\eta)=\eta$ and $\phi_q(\widetilde{m})=\widetilde{m}$,
and SBGD${}_{2,1}$ with $\phi_p(\eta)=\eta^2$ which seems to be a universally better.

\subsection{2D comparison with the Adam method}\label{sec:AdamandGD}
We  report on the results of SBGD${}_{pq}$, compared with  the Adam  method for the 2D Rastrigin
function \ref{eq:Rastrigin}. We also include  results for the \SGD and \MAGD methods. As in the case of the 1D objective function \eqref{eq:flat basins}, we distinguish between two cases:  initial data  uniformly distributed in $[-3,3]^2$ enclosing the global minimum at the origin, vs. initial data uniformly distributed at $[-3,-1]^2$.
The SBGD${}_{pq}$ variants with $(p,q)=(1,1)$ and $(p,q)=(2,1)$ were computed for $m=1000$ simulations with random with descent parameter $\lambda=0.8$, and using tolerance parameters, $tolmerge=0.1$,  $tolm=0.01$ and $tolres=10^{-4}$.\newline
 The SBGD${}_{pq}$  variant with $(p,q)=(2,1)$ seems consistently better than the vanilla version $(p,q)=(1,1)$. This will be further explored in the next \S\ref{sec:pq}.
When the methods are initiated at $[-3,3]^2$, the SBGD variants  provide comparable or  better results than  \SGD, \MAGD and Adams methods.
When initiated at $[-3,-1]^2$ which does not enclose $\bxmin$,  the SBGD variants provide distinctively better   results than  \SGD, \MAGD and Adam methods. The Adam iterations with the smaller initial step size $h_0=0.2$ remain trapped in local basins of attraction.

\begin{table}[h!]
\setstretch{1.5}
       \centering
    \begin{tabular}{|c||c|c|c|c|c|}
    \hline
    N&  5 & 10 & 15 & 20 & 30 \\
    \hline
           SBGD${}_{11}$ & 34.4\% & 52.1\% & 62.6\% & 70.0\% & 75.8\% \\
     \hline 
     SBGD${}_{21}$ & 34.5\% & 60.1\% & 75.3\% & 84.3\% & 91.0\% \\
     \hline 
       GD(0.004)  &  36.3\% & 50.5\% & 60.0\% & 70.0\% & 78.1\% \\
       \hline
       \MAGD  &  35.0\% & 51.0\% & 62.0\% & 70.8\% & 79.3\% \\
       \hline
       Adam(0.8)  &  23.7\% & 29.6\% & 39.1\% & 46.8\% & 65.5\% \\
       \hline
       Adam(0.2)  &  32.1\% & 40.9\% & 55.9\% & 65.3\% & 79.4\% \\
       \hline
    \end{tabular}
         \smallskip
 \caption{Success rates  of SBGD compared with \SGD,  \MAGD and Adam methods for global optimization of \eqref{eq:Rastrigin}($B=0$), based on   $m=1000$ runs with uniformly generated initial data  in $[-3, 3]^2$.}\label{tab:2D FB33 SBGD vs Adam GD}
 \end{table}

\begin{table}[h!]
\setstretch{1.5}
       \centering
    \begin{tabular}{|c||c|c|c|c|c|}
    \hline
    N&  5 & 10 & 15 & 20 & 30 \\
    \hline
      SBGD${}_{11}$ & 17.0\% & 49.2\% &  61.7\% & 67.0\% &72.7\% \\
     \hline 
     SBGD${}_{21}$ & 14.2\% & 46.7\% & 68.4\% & 81.9\% & 89.6\% \\
     \hline 
     GD(0.004) &  0.0\% & 0.0\% & 0.0\% & 0.0\% & 0.0\% \\
       \hline
       \MAGD  &  1.8\% & 2.4\% & 3.4\% & 4.3\% & 5.9\% \\
       \hline
              Adam(0.8)  &  24.5\% & 31.3\% & 41.4\% & 49.2\% & 66.9\% \\
       \hline
       Adam(0.2)  &  0.0\% & 0.0\% & 0.0\% & 0.0\% & 0.0\% \\
       \hline
    \end{tabular}
         \smallskip
 \caption{Success rates of  SBGD compared with \SGD, \MAGD  and Adam methods for global optimization of 2D Rastrigin 
  \eqref{eq:Rastrigin}($B=0$)  based on $m=1000$ runs of uniformly generated initial data in $[-3, -1]^2$.}\label{tab:2D FB31 SBGD vs Adam GD}
 \end{table}

\section{Numerical results --- $20$-dimensional problems}\label{sec:results-20D}
We now turn attention to  the computation of the global minimizer for the 20-dimensional Rastrigin and Ackley functions. As the dimension of the ambient space to be explored increases, so does the  number of agents, $N$, necessary  to explore that space in order to secure a `faithful'  approximate minimizer. The increase of  $N$ dependence on the dimension $d$ is intimately related to the way one quantifies the quality of such an approximation .

\subsection{Success rate --- $N$ vs. $d$} One approach to measure success that was used in one- and two-dimensional problems, is to secure a  computed solution within a pre-determined neighborhood of the global minimizer, $[\bxmin-0.25, \bxmin+0.25]^{d}$. This approach places severe restrictions in the case of high-dimensional data.
Specifically, Table \ref{tab:Ackley N of agents} and even more so, Table \ref{tab:Rastrigin N of agents}, show the rapid growth in the number of SBGD agents, $N=N(d)$, which are required  to ensure $80\%$ success rate in $m=500$ tests of high-dimensional Ackley, and respectively, 70\% success rate in Rastrigin benchmark functions.
 In both cases, one observes a rather small \emph{critical dimension}, $d_c$,  such that $N(d)\gg N(d_c)$ for $d > d_c$. 

\begin{table}[h]
\setstretch{1.5}
\vspace{0.1in}
    \begin{center}
    \parbox{.45\linewidth}
    {\begin{tabular}{|c||c|c|c|c|c|c|c|c|}
    \hline
        d & 10 & 11 & 12 & 13 & 14 & 15 & 16 \\
    \hline    
        N & 15 & 18 & 23 & 42 & 120 & 540 & 3000 \\
    \hline    
    \end{tabular}   
    \medskip
     \subcaption{80\% success rate with Ackley function}\label{tab:Ackley N of agents}}
     \parbox{.45\linewidth}
       {\begin{tabular}{c|c||c|c|c|c|}
    \cline{2-6}
        \hspace*{1.7cm} & d & 1 & 2 & 3 & 4  \\
    \cline{2-6}   
       \hspace*{1.7cm} & N & 4 & 23 & 180 & 2900 \\
    \cline{2-6}    
    \end{tabular}
        \medskip
      \subcaption{70\% success rate with Rastrigin function}}
      \caption{Success rate of $\bx_{SOL}\in [\bxmin-0.25,\bxmin+0.25]^{20}$.}\label{tab:Rastrigin N of agents}
      \end{center}
\end{table}

\begin{table}[h!]
\setstretch{1.3}
       \begin{center}
    \begin{tabular}{|m{1.5cm}|m{1.8cm}|| m{1.5cm} m{1.5cm} m{1.5cm}|}
    \hline
       $\bx_{*} = B$ &    &  N=50 & N=100 & N=200\\ [0.5ex]
       \hline
     $B = 0$  &SBGD${}_{21}$ \MAGD Adam(0.5)
       & 9.00e-07 1.18e-01 3.37e-03 & 2.02e-07  6.90e-02 4.03e-03 &      1.43e-07 1.88e-02 4.96e-03 \\ [0.5ex]
     \hline
     $B = 3$  &
     SBGD${}_{21}$ \MAGD Adam(0.5)  & 1.65e-06  7.64 2.14e-01 & 1.51e-06  6.72 1.22e-01&      7.96e-07 5.67 1.27e-01\\ [0.5ex]
    \hline          
         $B = 5$  &
    SBGD${}_{21}$ \MAGD Adam(0.5)  & 
     4.15  17.99 11.01 & 1.07  17.59 9.27 & 2.36e-01 17.11 7.52 \\   [0.5ex]  
    \hline          
    \end{tabular}
    \medskip
     \caption{$\EE[\bx_{SOL}-\bxmin\myr{]}$    for 20-dimensional  shifted Ackley. $m = 1000$.}\label{tab:Shifted Ackley error}
     \end{center}
\end{table}

\begin{table}[h!]
\hspace{0.05in}
\setstretch{1.3}
    \begin{center}
    \begin{tabular}{|m{1.5cm}|m{1.8cm}|| m{1.5cm} m{1.5cm} m{1.5cm}|}
    \hline
       $\bx_{*}=B$       &  & N=50 & N=100 & N=200\\ [0.5ex]
       \hline
     $B = 0$  &
     SBGD${}_{21}$ \MAGD Adam(0.5)  & 4.03e-01  2.96e-01 3.75e-01 & 4.98e-01 3.25e-01 2.86e-01 &      3.91e-01 4.35e-01 4.14e-01\\ [0.5ex]
     \hline
         $B = 3$ &
    SBGD${}_{21}$ \MAGD Adam(0.5)  &
     5.07   9.43 9.02 & 3.71  9.02  8.63 &      2.92 8.57 8.15\\ [0.5ex]
     \hline          
    $B = 5$ 
    &SBGD${}_{21}$ \MAGD Adam(0.5)  &
    3.53  18.02 17.13& 2.45  17.64 16.73 &      1.29 17.18 16.3\\ [0.5ex]
      \hline          
    \end{tabular}
    \medskip
        \caption{$\EE[\bx_{SOL}-\bxmin|$    for 20-dimensional  shifted Rastrigin. $m = 1000$.}\label{tab:Shifted Rastrigin error}
        \end{center}
\end{table}
\vspace{0.05cm}

An alternative approach to quantify the quality of  computed  minimizers is to measure the expected (average) error in position, $\EE[\bx_{SOL}-\bxmin|$. 
The results recorded in Tables \ref{tab:Shifted Ackley error} and \ref{tab:Shifted Rastrigin error} for 20-dimensional Ackley and, respectively, Rastrigin functions, indicate the advantage of SBGD${}_{21}$ as the shift $B$ increases. But the main point to observe is that even for relatively large $N=100,200$, the results  
fail to faithfully capture the global minimizer. Indeed, we claim that going beyond the two examples of Ackley and Rastrigin, measuring the (average) distance to of the computed minimizers, is not necessarily an effective quantifier for the quality of an optimizer: one might need a large $N=N(d)$ before approaching a small neighborhood of the  global minimizer.

\subsection{Measuring the loss} A more `faithful' way to measure the quality of numerical optimizers  in high-dimensional data is to measure the average loss (height), $\EE[F(\bx_{SOL}\myr{)}]$. After all, the underlying goal is to minimize the value of $F$. In particular, one might argue that $\bx_{SOL}$  will provide a faithful approximation whenever $F(\bx_{SOL})-F(\bxmin)$ is small, even  if $\bx_{SOL}$ remains far from $\bxin$. This approach of `looking from above', masks the increasing complexity with the increasing dimension.

\begin{table}[h!]
\setstretch{1.4}
\vspace{0.2in}
    \begin{center}
    \begin{tabular}{|m{1.5cm}|m{1.8cm}|| m{1.2cm} m{1.2cm} m{1.2cm} m{1.2cm} m{1.2cm} m{1.2cm}|}
    \hline
      $\bx^* = B$  & & N=5 & N=10  &  N=20& N=50& N=100 & N=200\\ [0.5ex]
        \hline
        $B = 0$  & SBGD${}_{21}$ \MAGD Adam(0.5)& 
         54.30 53.43 48.66 &  43.74 44.62 44.02 &  39.34 39.43 37.35 &  34.55 34.55 32.65 &  32.01 31.85 30.46 &  33.95 29.53 27.76\\ [0.5ex]
             \hline
        $B = 1$  &  SBGD${}_{21}$ \MAGD Adam(0.5)& 
        65.75 64.86 56.02 &  53.24 53.98 50.35 &  47.14 47.1 42.56 &   40.90 40.77 37.69 &   37.20 37.06 34.58 &  33.95 33.57 31.79\\ [0.5ex]
             \hline
        $B = 3$  &  SBGD${}_{21}$ \MAGD Adam(0.5)& 
        189.24 187.8 164.46 &  166.46 167.58 146.82 &  149.14 152.46 137.16 &   121.45 122.35 122.84 &   102.98 125.44 115.16 &  94.77 115.83 105.5\\ [0.5ex]
        \hline
        $B = 5$  &  SBGD${}_{21}$ \MAGD Adam(0.5)& 
       463.09 465.42 411.77 &  387.88 434.9 381.19 &  254.23 409.11 364.68 &   200.34 381.42 339.81 &   160.62 362.54 325.86 &  142.22 345.52  309.99\\ [0.5ex]
        \hline
    \end{tabular}
    \medskip
    \caption{20D Rastrigin. Average loss  with $m = 1000$ simulations.} \label{tab:2D Rastrigin height at $[-3,3]^2$}
    \end{center}
\end{table}

\begin{table}[h!]
\setstretch{1.4}
\vspace{0.2in}
    \begin{center}
    \begin{tabular}{|m{1.5cm}|m{1.8cm}||  m{1.5cm} m{1.5cm} m{1.5cm} m{1.5cm} m{1.5cm} |}
    \hline
      $\bx^* = B$  & &  N=10  &  N=20& N=50& N=100 & N=200\\ [0.5ex]
        \hline
        $B = 0$  & SBGD${}_{21}$ \MAGD Adam(0.5)& 
           0.89 2.76 0.37 &  2.62e-05 1.92 0.073 &  1.17e-05 0.97 0.061 &  9.72e-07 0.41 0.061 &  1.82e-07 7.73e-02 0.061\\ [0.5ex]
             \hline
        $B = 1$    & SBGD${}_{21}$ \MAGD Adam(0.5)& 
           2.38 3.01 0.87 &  0.27 2.00 0.17 &  3.9e-05 0.92 0.061 &  6e-06 0.31 0.061 &  1e-06 5,45e-02 0.061\\ [0.5ex]
             \hline
        $B = 3$    & SBGD${}_{21}$ \MAGD Adam(0.5)& 
           7.67 7.91 3.74 &  4.41 7.47 2.65 &  0.72 6.81 1.70 &  0.06 6.13 1.22 &  5.4e-05 5.37 0.79\\ [0.5ex]
        \hline
        $B = 5$    & SBGD${}_{21}$ \MAGD Adam(0.5)& 
           12.01 11.98 10.72 &  11.49 11.8 9.70 &  10.01 11.53 8.45 &  8.47 11.36 7.58 &  7.09 11.17 6.67\\ [0.5ex]
        \hline
    \end{tabular}
    \medskip
    \caption{20D Ackley. Average loss  with $m = 1000$ simulations.} \label{tab:2D Ackley height at $[-3,3]^2$}
    \end{center}
\end{table}

We begin, in Table \ref{tab:2D Rastrigin height at $[-3,3]^2$}, where  we record the average `loss'  obtained by SBGD${}_{21}$, \MAGD and Adam(0.5) with initial data equi-distributed in $[-3,3]^2$
for the 20-dimensional Rastrigin function. SBGD${}_{21}$ was used with descent parameter
$\clam=0.8$, shrinkage parameter $\gamma=0.5$ and thresholds
\[
tolmerge=\textnormal{1e-1}, \  \ tolm=\textnormal{1e-3}, \ \ tolres=\textnormal{1e-2}.
\]
The reason for the smaller shrinkage parameter $\gamma=0.5$ was  efficiency: backtracking is accelerated with more rough backtracking steps, yet  this does not seem to deteriorate the quality of  SBGD results.
As indicated before, the results for SBGD${}_{11}$ were only slightly worse than but otherwise comparable to SBGD${}_{21}$ and therefore are not recorded here. We make two observations.

\smallskip\noindent
(i) For  a small number of agents, $N\leq10$, the results,  particularly SBGD${}_{21}$ and \MAGD are comparable. Indeed, we recall that SBGD eliminates the lightest agents, so that after $N-1$ iterations, it is left with a single agent which explores the  large, uncharted ambient space, much like a single-agent method.\newline   
(ii) There is a clear trend that we saw before: SBGD${}_{21}$  outperforms the non-communicating \MAGD and Adam  when the global minimum is \emph{not} enclosed within the initial domain of initial data. While the results are comparable for $B=0,1$ and/or small $N$'s, there is an increasing  difference for $B=3,5$ and $N>20$.

\smallskip
Finally, in Table \ref{tab:2D Ackley height at $[-3,3]^2$} we report on the corresponding comparison of average loss for the 20-dimensional Ackley function. Here, one encounters a much more sensitive dependence on the initial step size: we had to increase $h_0=2$ (instead of $h_0=1$ used before) in order to realize the advantage of SBGD${}_{21}$. We maintain the usual backtracking parameters $\clam=0.2$ and $\gamma=0.9$, and threshold parameters
\[
tolmerge=\textnormal{1e-3}, \  \ tolm=\textnormal{1e-4}, \ \ tolres=\textnormal{1e-4}.
\]

\subsection{SBGD as pre-conditioner} A more practical strategy is to take the expected value, $\EE[\bx_{SOL}]$, varying overall SBGD solutions resulted from $m$ randomly generated initial configurations,  as a good initial guess and then iterate it with the steepest descent correction. In this way,  the algorithm is expected to mimic the convergence in expectation
$\displaystyle \EE[\bx_{SOL}] \stackrel{n\rightarrow \infty}{\longrightarrow} \bxmin$. Table \ref{tab:exp. errors} records the $L^{\infty}-$distance between the expectation $\EE[\bx_{SOL}]$ and the global minimizer $\bxmin$ for the  20-dimensional Ackley and Rastrigin functions. The expectation is computed with $m = 1000$ runs. We also present the error of the corrected solution $\bx_{corr}$, which is obtained by iterating $\EE[\bx_{SOL}]$ with gradient descent until $\displaystyle |\nabla F(\bx_{corr})|_{2}<10^{-3}$. For both benchmark functions, this strategy gives very good solutions with a reasonable amount of agents. The quality of $\EE[\bx_{SOL}]$ is improved with more agents applied.

\vspace{0.3in}
\begin{table}[h]
\setstretch{1.5}
    \centering
    \begin{tabular}{|m{3.3cm}|m{3cm}|| m{1.6cm} m{1.6cm} m{1.8cm}|}
    \hline
                   &  & N=50 & N=100 & N=200\\ [0.5ex]
     \hline
    Ackley  function  & 
    $|\EE[\bx_{SOL}]-\bx_{*}|_{\infty}$ $|\bx_{corr}-\bx_{*}|_{\infty}$  &  
     	$8.29e^{-2}$ $2.07e^{-13}$ & $7.28e^{-2}$ $1.96e^{-13}$ & 
     	$6.79e^{-2}$ $2.18e^{-13}$\\ [0.5ex]
    \hline          
    Rastrigin  function  &
    $|\EE[\bx_{SOL}]-\bx_{*}|_{\infty}$ $|\bx_{corr}-\bx_{*}|_{\infty}$ &
     $2.07e^{-1}$ $9.06e^{-6}$ & $7.86e^{-2}$ $2.29e^{-6}$ &
     $7.36e^{-2}$  $2.05e^{-6}$ \\    [0.5ex]
    \hline          
    \end{tabular}  
    \medskip
    \caption{Errors of the expectation and the corrected solution in $d = 20$, $B = C = 0$.}\label{tab:exp. errors}
\end{table}

We also investigate the effect of shifting in the initial distribution. The solution $\bx_{corr}$ is considered to be a correct approximation of global minimizer if it lies in the region $Q_{0.5}(\bx_{*})$. Table \ref{tab:Shifted Ackley error} and \ref{tab:Shifted Rastrigin error} show the results for the Ackley and the Rastrigin function under varying shift parameters and different agent numbers. It turns out that the quality of the solution is very sensitive to the initial data in high dimension. The algorithm fails to work correctly in the Ackley test cases when $B = 2$. The situation is even worse for the Rastrigin test cases due to the unclear difference between different minimums. The method gives wrong solutions in all the test cases when $B$ is 1.5.

\begin{table}[h!]
\setstretch{1.5}
       \begin{center}
    \begin{tabular}{|m{1.5cm}|m{3cm}|| m{1.6cm} m{1.6cm} m{1.8cm}|}
    \hline
       $\bx_{*} = B$       &  &  N=50 & N=100 & N=200\\ [0.5ex]
     \hline
     $B = 1$  &
     $|\EE[\bx_{SOL}]-\bx_{*}|_{\infty}$ $|\bx_{corr}-\bx_{*}|_{\infty}$ & $3.97e^{-1}$ $2.10e^{-13}$ & $3.53e^{-1}$ $1.29e^{-13}$ & 
     $3.36e^{-2}$ $1.15e^{-13}$\\ [0.5ex]
    \hline          
     $B = 1.5$  &
    $|\EE[\bx_{SOL}]-\bx_{*}|_{\infty}$ $|\bx_{corr}-\bx_{*}|_{\infty}$  & 
    	$7.45e^{-1}$ $9.69e^{-1}$ & $6.78e^{-1}$ $1.27e^{-13}$ &
    	$5.43e^{-1}$ $1.60e^{-13}$  \\   [0.5ex]
    \hline          
    $B = 2$  &
    $|\EE[\bx_{SOL}]-\bx_{*}|_{\infty}$ $|\bx_{corr}-\bx_{*}|_{\infty}$ &
     $1.20$ $9.69e^{-1}$ & $1.09$ $9.69e^{-1}$ &$1.01$ $9.68e^{-1}$\\   [0.5ex]  
    \hline          
    \end{tabular}
    \medskip
     \caption{Solutions for the shifted Ackley function by SBGD, $d = 20$, $m = 1000$.}\label{tab:Shifted Ackley error}
     \end{center}
\end{table}

\begin{table}[h!]
\hspace{0.2in}
\setstretch{1.5}
    \begin{center}
    \begin{tabular}{|m{1.5cm}|m{3cm}|| m{1.6cm} m{1.6cm} m{1.8cm}|}
    \hline
       $\bx_{*}=B$       &  & N=50 & N=100 & N=200\\ [0.5ex]
     \hline
     $B = 0.5$  &
     $|\EE[\bx_{SOL}]-\bx_{*}|_{\infty}$ $|\bx_{corr}-\bx_{*}|_{\infty}$ &
      $2.21e^{-1}$  $5.86e^{-6}$ &$2.38e^{-1}$ $1.06e^{-5}$ & 
     $2.24e^{-1}$ $6.53e^{-6}$ \\ [0.5ex]
     \hline          
    $B = 1$ &
    $|\EE[\bx_{SOL}]-\bx_{*}|_{\infty}$ $|\bx_{corr}-\bx_{*}|_{\infty}$ &
     $5.97e^{-1}$ $1.99$  & $4.80e^{-1}$ $5.89e^{-6}$  &
     $4.19e^{-1}$ $1.11e^{-5}$\\ [0.5ex]
     \hline          
    $B = 1.5$ 
    &$|\EE[\bx_{SOL}]-\bx_{*}|_{\infty}$ $|\bx_{corr}-\bx_{*}|_{\infty}$ &
     $7.92e^{-1}$ $9.95e^{-1}$ & $7.00e^{-1}$ $9.95e^{-1}$ &
     $6.34e^{-1}$ {$1.99$}\\ [0.5ex]
      \hline          
    \end{tabular}
    \medskip
        \caption{Solutions for the shifted Rastrigin function by SGBD, $d = 20$, $m = 1000$.}\label{tab:Shifted Rastrigin error}
        \end{center}
\end{table}
\vspace{0.2in}

The results above indicate that the solutions obtained for high-dimensional problems can be significantly affected by the shifting in the target function, especially when the global minimum is not very distinguishable from the other minima. 

\bibliographystyle{plain}
\bibliography{Ref}

\end{document}